\documentclass[11pt]{article}
\usepackage{amsmath,epsfig,amssymb,amsfonts}
\usepackage[bottom,PostScript=dvips]{diagrams}

\newcounter{SEC}

\numberwithin{equation}{section}

\newcounter{STAT}[section]
\renewcommand{\theSTAT}{\arabic{section}.\arabic{STAT}}
\renewcommand{\subsection}[1]{\refstepcounter{STAT}%
\par\vspace{2em plus .2em minus .1em}\par%
\noindent {\large\bfseries \theSTAT. \ #1. \ }%
}

\newenvironment{Thm}{\refstepcounter{STAT}%
\par\vspace{1.5ex}\par%
\noindent {\bf Theorem \theSTAT.}\begingroup\itshape\  }{ \endgroup\par\vspace{.8ex}\par}
\newenvironment{Defn}{\refstepcounter{STAT}%
\par\vspace{1.5ex}\par%
\noindent {\bf Definition \theSTAT.}\begingroup\itshape\  }{ \endgroup\par\vspace{.8ex}\par}
\newenvironment{Cor}{\refstepcounter{STAT}%
\par\vspace{1.5ex}\par%
\noindent {\bf Corollary \theSTAT.}\begingroup\itshape\  }{ \endgroup\par\vspace{.8ex}\par}
\newenvironment{Lem}{\refstepcounter{STAT}%
\par\vspace{1.5ex}\par%
\noindent {\bf Lemma \theSTAT.}\begingroup\itshape\  }{ \endgroup\par\vspace{.8ex}\par}
\newenvironment{Prop}{\refstepcounter{STAT}%
\par\vspace{1.5ex}\par%
\noindent {\bf Proposition \theSTAT.}\begingroup\itshape\  }{ \endgroup\par\vspace{.8ex}\par}
\newenvironment{Rem}{\refstepcounter{STAT}%
\par\vspace{1.5ex}\par%
\noindent {\bf Remark \theSTAT.}\begingroup\itshape\  }{ \endgroup\par\vspace{.8ex}\par}
\newcommand{\set}[1]{\left\{#1\right\}}

\newcommand{\CC}{\mathbb C}
\newcommand{\NN}{\mathbb N}
\newcommand{\ZZ}{\mathbb Z}
\newcommand{\PP}{\mathbb P}
\newcommand{\QQ}{\mathbb Q}
\newcommand{\m}{\mathfrak M}
\newcommand{\eps}{\varepsilon}
\newcommand{\vf}{\varphi}
\newcommand{\de}{\delta}
\newcommand{\To}{\longrightarrow}
\newcommand{\xTo}[1]{\xrightarrow{#1}}
\newcommand{\A}{\mathcal{A}}
\newcommand{\T}{\mathcal{T}}
\newcommand{\M}{\mathcal{M}}
\renewcommand{\O}{\mathcal{O}}

\newcommand\E{{\mathcal E}}
\newcommand{\Vect}{\mathcal Vect}

\newenvironment{proof}[1]{\par\smallskip\par\noindent{\bf
Proof#1}}{$\ \scriptstyle \blacksquare$\par\bigskip\par}

\newcommand\arcctg{\mathrm{arcctg}\,}
\renewcommand\Im{\mathop\mathrm{Im}}

\newcommand{\Ker}{\mathop\mathrm{Ker}}
\newcommand{\Coker}{\mathop\mathrm{Coker}}
\newcommand{\Coh}{\mathcal{C}oh\,}

\newcommand{\Hom}{\mathrm{Hom}\,}
\newcommand{\Ext}{\mathrm{Ext}\,}

\newcommand{\x}{\times }
\newcommand\Aut{\mathop\mathrm{Aut}}

\newcommand\rk{\mathop\mathrm{rk}}
\newcommand{\seq}[3]{0\To #1\To #2\To #3\To 0}
\renewcommand{\le}{\leqslant}
\renewcommand{\ge}{\geqslant}

\newcommand{\TR}[3]{%
\TRM{#1}{#2}{#3}{}{}}

\newcommand\TRM[5]{\begin{diagram}[size=1.5em]
&&#3\\
&\ruTo^{\scriptstyle #5}&&\rdDashto\\
#2&&\lTo_{\scriptstyle #4}&&#1
\end{diagram}}

\newcommand\TFILTR[7]{%
\begin{diagram}[size=1.4em]
                         &                           &{#2}_{{#5}_0}             &         &          &                           &{#2}_{{#5}_1}             &         &          &          &      &    &&&{#2}_{{#5}_{#4}}\\
                         &\ruTo^{\scriptstyle {#6}_0}&                          &\rdDashto&          &\ruTo^{\scriptstyle {#6}_1}&                          &\rdDashto&          &          &      &    &&\ruTo^{\scriptstyle {#6}_{#4}}&&\rdDashto\\
\llap{${#1}=$}{#3}^0{#1} &                           &\lTo_{\scriptstyle {#7}_1}& &{#3}^1{#1}&
&\lTo_{\scriptstyle {#7}_2}&
&{#3}^2{#1}&\lTo\relax&\cdots&\lTo&{#3}^{#4}{#1}&&\lTo_{\scriptstyle
{#7}_{{#4}+1}}&&{#3}^{{#4}+1}{#1}=0
\end{diagram}}

\newcommand\TFILT[5]{%
\begin{diagram}[size=1.4em]
&&{#2}_{{#5}_0}&&&&{#2}_{{#5}_1}&&&&&&&&{#2}_{{#5}_{#4}}\\
&\ruTo&&\rdDashto&&\ruTo&&\rdDashto&&&&&&
\ruTo&&\rdDashto\\
\llap{${#1}=$}{#3}^0{#1}&&\lTo&&{#3}^1{#1}&&\lTo&&
{#3}^2{#1}&\lTo\relax&\cdots&\lTo&{#3}^{#4}{#1}&&\lTo&&{#3}^{{#4}+1}{#1}\rlap{$=0$}
\end{diagram}}

\newcommand\TF[5]{%
\begin{diagram}[size=1.4em]
                         &                           &{#2}_{_0}             &         &          &                           &{#2}_{1}             &         &          &          &      &    &&&{#2}_{{#4}}\\
                         &\ruTo^{\scriptstyle {#5}_0}&                          &\rdDashto&          &\ruTo^{\scriptstyle {#5}_1}&                          &\rdDashto&          &          &      &    &&\ruTo^{\scriptstyle {#5}_{#4}}&&\rdDashto\\
\llap{${#1}=$}{#3}^0{#1} &                           &\lTo & &{#3}^1{#1}& &\lTo &
&{#3}^2{#1}&\lTo\relax&\cdots&\lTo&{#3}^{#4}{#1}&&\lTo &&{#3}^{{#4}+1}{#1}=0
\end{diagram}}

\newcommand{\LA}{\left\langle}
\newcommand{\RA}{\right\rangle}

\newcommand{\DST}[3]{\left(#1,\{{#3}_{#2}\}_{#2\in #1}\right)}
\newcommand{\DSTm}[3]{\left(#1,\{{#3}_{#2}\}_{#2\in #1},\tau_{#1}\right)}

\newcommand{\PGN}{HN-system}
\newcommand{\HN}{HN-system}

\newcommand{\tr}{distinguished triangle}

\newcommand{\com}{{\scriptscriptstyle\bullet}}
\begin{document}

\begin{flushright}
%\itshape In remembrance  of\\
%\bfseries A.\,N. Tyurin\rule{18pt}{0pt}
{\small\itshape\bfseries To blessed memory of  Andrey Tyurin\/}
\end{flushright}

\vspace*{3mm}

\begin{center}
{\LARGE\bfseries Stability data and t-structures on\\ a triangulated category}\\[5mm]
{\large\itshape A.~Gorodentscev, S.~Kuleshov\footnote{The second author was
    partially supported by INTAS  grant OPEN 2000 269}, A.~Rudakov}
\end{center}

\vspace*{1mm}

\begin{center}\parbox{100mm}{\small \hspace*{1em}
 We propose the notion of stability on  a
triangulated category that is a generalization of the   T.~Bridgeland's stability data. We
establish connections between stabilities and t-structures on a category and as application we
get the classification of bounded t-structures on the category $D^b(\Coh\PP^1)$.}\end{center}

\vspace*{3mm}

% --

\section{Introduction}

\noindent  For quite a number of years the authors believed that certain generalization of
stability would be desirable to use in the context of derived and triangulated categories, but
could not come to the satisfactory definition.

Recently in his article \cite{BR1}, T.Bridgeland, following not so rigorously presented, but
inspiring physical ideas of Michael R. Douglas (\cite{D}) provided the definition of stability
for a triangulated category. Bridgeland also showed how the key properties of stability can be
reformulated in this context and worked on describing all the stabilities for a given category
via a kind of moduli space of stabilities. Although it is not clear to us where the approach to
constructing the moduli space of stabilities proposed by T.Bridgeland will lead, we believe that
his definition of the stability makes the breakthough.

In this paper we propose the definition of stability for a triangulated category or, in short,
t-stability which is the modification and in fact the generalization of the definition given by
T.\,Bridgeland. We believe we have taken away ``unnecessary details'' keeping the core features
intact. In short, we exclude all about the ``central charge'' (in the sense of \cite{BR1}) from
the definition of t-stability, and we do not demand that the semi-stable subcategories are
ordered according the real numbers assigned to them as indexes. But we keep the way to
generalize the Harder-Narasimhan filtration that was proposed by T.Bridgeland.

We begin in Section \ref{Abelian category} reminding the basic properties of the stability for
an abelian category. Section \ref{DEFINITION} is devoted to the definition of t-stability for a
triangulated category. We also discuss basic properties of t-stability and several basic
examples or constructions of a t-stability on a triangulated (or derived) category.

We show (p. \pageref{Gieseker stability to t-stability})  that the natural Gieseker  stability
for torsion free coherent sheaves extend to the derived category of coherent sheaves, and that
the helix generalization (\cite{GK}) of the Beilinson theorem about the basis of the derived
category $D^b(\Coh \PP^n)$ of coherent sheaves on the projective space $\PP^n$ (\cite{BEI})
leads to the remarkable t-stabilities for this category (Proposition \ref{Helstab}).

In Section \ref{PROPERTIES OF PGN} we discuss the Postnikov systems that we start to call in our
context by \mbox{t-fil\-t\-rations}. We show that the choice of t-stability provides each object
with the canonical Harder-Narasimhan t-filtration   and we  describe properties of these
HN-filtrations.

In Section \ref{ORDER} we develop  several methods to connect t-stabilities and t-structures on
a category, as well as to make more coarse or refined t-stabilities.

For the case $\PP^1$ this enables us to get the full classification of t-stabilities on
$D^b(\Coh \PP^1)$ that we present in Subsection \ref{T-STABILIIES FOR P1}. As a consequence in
Subsection \ref{CLASS T-STRUCT P1} we obtain the full list of bounded  t-structures on $D^b(\Coh
\PP^1)$ (considered as a triangulated category).

In Section \ref{ELLIPTIC} we show how to describe all t-stabilities for coherent sheaves on an
elliptic curve.

\section{Stability data on abelian categories}\label{Abelian category}

\noindent  We begin with a remark about the  stability   of coherent sheaves. This notion arose
as a tool for construction  moduli spaces of coherent torsion free sheaves on varieties and came
from the invariant theory. The  moduli space is obtained as an orbit space for a reductive group
action on a vector space that parameterizes   sheaves with fixed topological invariants. In this
approach moduli space points corresponding to closed orbits are well defined. Closed orbits and
sheaves containing in such orbits are called stable. Unfortunately, in many cases the set of
stable (closed) orbits   is not compact. To compactify it we should add   orbits, whose closure
  does not contain the zero vector. Such orbits and corresponding sheaves are called
semistable.

 There exist  numerical criteria  of stability for orbits (and sheaves) (see, for example
 \cite{GIE, MA1}, and \cite{MA2}).
Roughly speaking, we correspond a vector  $\mu(F)$ (or, simply, number) to each a torsion free
sheaf $F$. This vector is called a   slope  of the sheaf. The criterium says that a torsion free
sheaf $F$ is semistable iff for any nonzero subsheaf $E\subset F$ we have $\mu(E)\le\mu(F)$ (the
vectors are compared  lexicographically). The stability of a sheaf $F$ means that $F$ is
semistable and simple (i.e. $\Hom(F,F)=\CC$).

As a rule, each a component of a vector  slope is a ratio of additive functions on $K_0(\Coh)$,
the Grothendieck group  of a coherent sheaves category. Therefore slope satisfies the following
{\itshape seesaw condition\/}:
\begin{quote}\itshape
for an exact sequence of torsion free sheaves
\[\seq EFG\]
we have
\[\begin{split}
\mu(E)<\mu(F)\quad\Leftrightarrow\quad\mu(F)<\mu(G),\label{stability define properties}\\
\mu(E)=\mu(F)\quad\Leftrightarrow\quad\mu(F)=\mu(G),\\
\mu(E)>\mu(F)\quad\Leftrightarrow\quad\mu(F)>\mu(G).
\end{split}\]
\end{quote}
Besides the moduli space construction, stability of coherent sheaves have two more useful
applications, following from the seesaw condition:
\begin{quote}\itshape
(i) for any semistable sheaves $E$ and $F$ an inequality $\mu(E)>\mu(F)$ implies\\
\hspace*{1.5em}\mbox{$\Hom(E,F)=0$};

(ii) any torsion free sheaf $X$ has a canonical Harder--Narasimhan filtration
\[
\begin{diagram}[size=1.5em]
\llap{$X=$}F^0X&\lInto&F^1X&\lInto&F^2X&\lInto\relax&\quad\cdots\quad&\lInto&F^nX&\lInto&F^{n+1}X\rlap{$=0$\ ,}\\
\dOnto&&\dOnto&&\dOnto&&&&\dOnto&&\\
G_1&&G_2&&G_3&&&&G_n
\end{diagram}
\]
where each vertical epimorphism  is a part of  the short exact sequence
\[\seq{F^{i}X}{F^{i-1}X}{G_i}\] with semistable $G_i$  and $\mu(G_i)<\mu(G_j)$ whenever $i<j$.
\end{quote}

So, we see that a stability of coherent sheaves gives a powerful filtration  but for  torsion
free sheaves only. We would like to define a stability in such a way that any nonzero object has
a Harder--Narasimhan filtration. The first abstract definition of stability on an abelian
category    was  done in \cite{RU} as follows.

\begin{Defn}\label{order-stability}
Let us say that a stability structure on an abelian category $\A$ is given if there is a
preorder on $\A$ such that for an exact sequence of nonzero objects $\seq A B C$ we have the
seesaw property:
\begin{align*}
 &\text{either}&&A\prec B\ \Leftrightarrow\ A\prec C\ \Leftrightarrow\ B\prec C,\\
 &\text{or}&&A\succ B\ \Leftrightarrow\ A\succ  C\ \Leftrightarrow\ B\succ   C,\\
 &\text{or}&&A\asymp B\ \Leftrightarrow\ A\asymp C\ \Leftrightarrow\ B\asymp C.
\end{align*}
We'll call a nonzero object $A\in\A$ semistable if $B\preceq A$ whenever  $0\ne B\subset A$.
\end{Defn}

In   \cite{RU} the following theorem was proved.

\begin{Thm}\label{Rudakov} 1. If objects $A$, $B$ are semistable and $A\prec B$, then
$\Hom(B,A)=0$.

2. Suppose $\A$ is weakly-artinian and weakly-noetherian, then for an object $X\in\A$ there
exists a unique Harder--Narasimhan filtration
\[\begin{diagram}[size=1.5em]
\llap{$X=$}F^0X&\lInto&F^1X&\lInto&F^2X&\lInto\relax&\quad\cdots\quad&\lInto&F^nX&\lInto&F^{n+1}X\rlap{$=0$}\\
\dOnto&&\dOnto&&\dOnto&&&&\dOnto&&\\
G_1&&G_2&&G_3&&&&G_n
\end{diagram}\]
with semistable quotients $G_i$'s such that $G_i\prec G_j$ whenever $i<j$.
\end{Thm}

Thus, having a preorder on an abelian category  satisfying the seesaw property and finiteness
conditions, we obtain the set of semistable objects and Harder--Narasimhan filtration for each
object. But there appears a question: how one can order the objects? In all known examples such
an order is obtained with the help of slope. Therefore  we propose an abstract definition of
slope on an abelian category, generalizing Bridgeland's central charge.

\begin{Defn}\label{abstract slope}  Let  $\A$ be an abelian category.  A  linearly independent  system of
additive functions $(x_0,\ldots, x_{r-1})$ on $K_0(\A)$ (the Grothendieck group   of $\A$) is
called positive   if for  any   $A\in\A$ the following   holds
\begin{equation}
 \begin{split} & x_0(A)\ge 0,\qquad\text{and}\\
  &x_0(A)=0  \Rightarrow x_1(A)\ge 0,\qquad\text{and}\\
  &x_0(A)=x_1(A)=0  \Rightarrow x_2(A)\ge 0,\qquad\text{and}\\
  &\dots\\
  &x_0(A)=\cdots=x_{r-2}(A)=0  \Rightarrow x_{r-1}(A)> 0.\\
   \end{split}
\end{equation}
If, in addition \[x_0(A)=\cdots=x_{r-1}(A)=0\  \Rightarrow\ A=0,\] the positive system is called
a positive base.

 Let $s=\min\limits_{i}\set{x_i(A)\ne 0}$ and
\[\gamma(A)=\left(\underbrace{1,\ldots,1}_s,\nu\left(-\frac{x_{s+1}(A)}{x_s(A)}\right),
\nu \left(-\frac{x_{s+2}(A)}{x_s(A)}\right),\ldots, \nu
\left(-\frac{x_{r-1}(A)}{x_{s}(A)}\right)\right),\] where
\[ \nu\left(\frac{a}{b}\right)=
  \begin{cases}
    \frac 1 \pi \arcctg\frac a b, & b>0, \\
    1 & b=0.
  \end{cases}
\] We call  $\gamma(A)$  the slope of an object $A\in\A$, determined by the positive system
 $(x_0,\ldots, x_{r-1})$ (not necessary a base).
\end{Defn}

For example, the base $(\rk, \deg)$  on the category $\Coh C$ of coherent sheaves on a smooth
algebraic curve $C$ has the positivity property, the base $(\rk, \deg, \chi(\O_S,\cdot))$  on
the category $\Coh S$ of coherent sheaves on a smooth algebraic surface $S$ with Picard's number
1 has the positivity property as well.

Since the slope $\gamma$ is formed via additive functions on $K_0(\A)$, the   ordering induced
by $\gamma$ ($A\preceq B$ $\Leftrightarrow$ $\gamma(A)\le \gamma(B)$) satisfies the seesaw
property. Therefore, we obtain a stability structure and, consequently, Harder--Narasimhan
filtration for each   object.

Another way to consider  stability on an abelian category is to take the properties (i) and (ii)
above as a definition. Namely,

\begin{Defn}\label{def stab on A} Let $\A$ be an abelian category and $\Phi$ be a linearly
 ordered  set. Suppose that for
any $\vf\in\Phi$ a  subcategory $\Pi_\vf\subset\A$ is determined and $\Pi_\vf$ are closed under
extensions. If the following properties are valid
\begin{enumerate}
  \item[(i)] $\Hom_{\A}(\Pi_{\vf'},\Pi_{\vf''})=0$ for $\vf'>\vf''$;
  \item[(ii)] each nonzero object $X\in\A$ has a Harder--Narasimhan filtration
  \[\begin{diagram}[size=1.5em]
\llap{$X=$}F^0X&\lInto&F^1X&\lInto&F^2X&\lInto\relax&\quad\cdots\quad&\lInto&F^nX&\lInto&F^{n+1}X\rlap{$=0$}\\
\dOnto&&\dOnto&&\dOnto&&&&\dOnto&&\\
G_1&&G_2&&G_3&&&&G_n
\end{diagram}\] with
  $G_i=F^iX/F^{i+1}X\in \Pi_{\vf_i}$ and $\vf_i<\vf_j$ for $i<j$;
\end{enumerate}
then we call the data $\DST \Phi \vf \Pi$   stability data (or stability) on $\A$.
\end{Defn}

Thus we come to three  definitions of stability on an abelian category:
\begin{enumerate}
  \item via ordering objects (let us call it order-stability);
  \item via base of $K^*(\A)$ with positivity property  (slope-stability);
  \item via collection of semistable categories (stability data).
\end{enumerate}

It is obvious, that any slope-stability on an abelian category induces an order-stability, and
the last one induces a stability data. But   we don't know: if one  can   reconstruct the
slope-stability (or the order-stability) starting with   given stability data on an abelian
category. The question is   interesting.

%Whether it is possible to    можно ли

\section{Definition and basic examples}\label{DEFINITION}

\noindent In Section \ref{Abelian category} we considered a   stability on an abelian
categories. Keeping in mind the three definitions, we see that the best for  extending  a
stability on a triangulated category is the last one (def. \ref{def stab on A}).
 We need the  generalization of the Harder--Narasimhan filtration to a triangulated category.
The first step in this direction was done by T.~Bridgeland in \cite{BR1}. We modify the
Bridgeland's definition of stability data, and exclude any reference to a central charge in
order to make them more general.

The classical  notion of a filtration  is based on subobjects and quotients. For a triangulated
category  Bridgeland proposed a natural generalization (\cite{BR1}), namely   a Postnikov's
system
\[\TFILTR X X F n {\vf} {q} {p}\]
(where each a triangle $\TRM{F^{i+1}X}{F^iX}{X_{\vf_i}}{p}q$\qquad is distinguished).  We call
such a system  a filtration of an object $X$ in  a triangulated category or in short "{\itshape
t-filtration\/}". It is natural to call the objects $X_i$  quotients and $F^iX$ terms of the
t-filtration.  Often we need information only about quotients of a t-filtration on $X$. In this
case we write: \quad $X\rightsquigarrow(X_0,X_1,\ldots, X_n)$ \quad as a notation for the
t-filtration.

We are ready to formulate the main definition.

\begin{Defn}\label{st-data}
Let   $\T$ be a triangulated category, $\Phi$ be a linearly ordered set.  Suppose  that for each
$\vf\in\Phi$ a strongly full, extension-closed\footnote{An extension-closed subcategory
$\A\subset\T$ means  that whenever $A\To B\To C\To A[1]$ is a \tr, with $A\in\A$ and $C\in \A$,
then $B\in\A$ also.} nonempty subcategory $\Pi_\vf\subset \T$ is defined. The pair
$\left(\Phi,\{\Pi_\vf\}_{\vf\in\Phi}\right)$ is called   {\itshape stability data\/} (or simply
{\itshape t-stability\/}) if
\begin{enumerate}
\item the shift of the triangulated category acts  on the set $\{\Pi_\vf\}_{\vf\in\Phi}$
in the following sense: there exists   $\tau\in{\rm Aut}\,\Phi$ such that
$\Pi_\vf[1]=\Pi_{\tau(\vf)}$ and $\tau(\vf)\ge\vf$;
\item   $\forall\,\psi>\vf\in \Phi$ \quad $\Hom^{\le
0}(\Pi_\psi,\Pi_\vf)=0 $;
\item for any nonzero object $X\in\T$ there exists a finite Postnikov's system:
\[\TFILTR X X F n {\vf} {q} {p}\]
with nonzero $X_{\vf_i}\in\Pi_{\vf_i}$ and $\vf_i<\vf_{i+1}$.
\end{enumerate}
We shall call  such a system  Harder--Narasimhan system (\PGN\ or  HN-filtration) of $X$, the
objects $X_{\vf_i}$ are $\Phi$-semistable quotients   (w.r.t stability data
$\left(\Phi,\{\Pi_\vf\}_{\vf\in\Phi}\right)$), the subcategories $\Pi_\vf$ are semistable
subcategories of slope $\vf$.
\end{Defn}

Fairly often we   denote   stability data $\DST \Phi \vf \Pi$ simply by $\Phi$.

 To construct
some trivial, but important examples, recall  (\cite{GM}) that a t-struc\-tu\-re on a
triangulated category $\T$ is a pair $\left( D^{\le0},D^{\ge0}\right)$ of strictly  full
subcategories, satisfying the following conditions:

\begin{enumerate}
\item $D^{\le 0}\subset D^{\le 0}[-1]$ and $D^{\ge 0}\supset D^{\ge 0}[-1]$;
\item $\Hom^0(D^{\le 0},D^{\ge 0}[-1])=0$;
\item $\forall\,X\in \T$ \  there exists a \tr\
$$\TRM{X_{\le0}}{X}{X_{\ge1}}pq$$
 with $X_{\ge1}\in D^{\ge0}[-1]$ and $X_{\le0}\in D^{\le
0}$.
\end{enumerate}

  If the following property holds, then the t-structure is called {\itshape bounded\/} (\cite{BR1}):
\begin{enumerate}
\item[4)] $\forall\,X\in \T$ $\exists\, m,\,n\in\ZZ$ such that  $X\in D^{\ge 0}[-m]\cap D^{\le
0}[-n]$
\end{enumerate}

  We shall  use the standard  notation: $D^{\le n}=D^{\le 0}[-n]$ and $D^{\ge
n}=D^{\ge 0}[-n]$.

\begin{Lem}\label{trivDST} Suppose $\left( D^{\le0},D^{\ge0}\right)$
is a t-structure on $\T$  such that   $D^{\le 0}= D^{\le 1}=D^{\le 0}[-1]$, $D^{\ge 0}= D^{\ge
1}= D^{\ge 0}[-1]$. Let
\[ \Phi=\{0,1\},\qquad
    \Pi_0=D^{\ge0}[-1],\qquad\Pi_1=D^{\le0},\qquad \tau(0)=0,\quad\tau(1)=1.
\]
Then $\DST \Phi \vf \Pi$ determines  stability data on $\T$.
\end{Lem}
\begin{proof}{} follows immediately from the definitions.\end{proof}

\begin{Lem}\label{stability by t-str} Let $(D^{\le0},D^{\ge0})$ be
a bounded  t-structure on a triangulated category $\T$.\\ Let \ $\Pi_i=\A[i]=D^{\le-i}\cap
D^{\ge-i}$, then $\left(\ZZ,\{\Pi_i\}_{i\in\ZZ}\right)$ makes  stability data on $\T$.
\end{Lem}
\begin{proof}{} can go as follows:
\begin{enumerate}
  \item[] Notice that  if $p<q$, then
  $\Hom^{\le0}(D^{\le-p},D^{\ge-q})=0$.
  \item[] Since the t-structure   $(D^{\le0},D^{\ge0})$ is bounded, then for each nonzero object
  $X\subset \T$ there exist  $n_+(X)\ge n_-(X)\in\ZZ$ such that
\begin{align*}
  &\Hom^0(X,D^{\ge-n_-(X)})\ne0, \qquad \Hom^{\le0}(X,D^{\ge-n})=0\text{ for } n<n_-(X); \\
  &\Hom^0(D^{\le-n_+(X)},X)\ne 0,  \qquad  \Hom^0(D^{\le-k},X)= 0\text{ for } k>n_+(X).
\end{align*}
  \item[] Denote  $n_-(X)$ by $n_0$ and consider a \tr\
\[
    \TR{Y_0}X{X_0}
\]
with $X_0\in D^{\ge-n_0}$, $Y_0\in D^{\le -(n_0+1)}$ (it exists by the definition of
t-structure). Clearly we get
\begin{itemize}
  \item $\Hom^{\le0}(Y_0,X_0)=0$;
  \item $X_0\ne0$;
  \item $X_0\in D^{\ge-n_0}\cap D^{\le-n_0}=\Pi_{n_0}$,
\end{itemize}
and $Y_0$ serves as the first term of the Harder--Narasimhan filtration.
\end{enumerate}
Now one can finish the proof   by induction.\end{proof}

We would like to note that  Definition \ref{def stab on A} and Lemma \ref{stability by t-str}
allow us to extend a stability from an abelian category $\A$ onto the bounded  derived category
of $\A$. Let us formulate slightly  more general  fact.

\begin{Prop}\label{t-stab from core}
 Let $\T$ be a triangulated category with a bounded t-structure $(D^{\ge0},D^{\le 0})$
 ({t-category}). Suppose that on the core $\A=D^{\ge0}\cap D^{\le 0}$ of the t-structure
  stability data $\DST \Phi \vf \Pi$
 are given. Consider the set $\Psi=\ZZ\times \Phi$ with the lexicographic order by $\Psi$
 put $P_{(i,\vf)}=\Pi_\vf[i]$. Then $\DST \Psi {(i,\vf)} P$ constitutes
  stability data on the category  $\T$.
\end{Prop}
\begin{proof}{.} There is only one non-obvious moment in the proof. Namely, the existence of the
finite \HN\ for each nonzero object. But this   follows directly  from Definition \ref{def stab
on A} of   stability data on an abelian category, Lemma \ref{stability by t-str} and Proposition
\ref{atom. and glu. t-filtr} that we  prove in the next section.%
\end{proof}

%%%%%%%%%%%%%%%%%%%%%%%

\begin{Cor}\label{Gieseker stability to t-stability}
Any stability structure $\DST \Phi \vf \Pi$ on an abelian category $\A$ induces a t-stability
$\left(\ZZ\x\Phi,\{\Pi_\vf[i]\}_{(i,\vf)\in\ZZ\x\Phi} \right)$ on the bounded derived category
$D^b(\A)$.
\end{Cor}

It was be shown in article \cite{RU} that Gieseker stability of torsion free coherent sheaves on
a projective variety $X$ one can extend to a stability structure on the category $\Coh X$ of all
coherent sheaves on $X$. Therefore due to the previous corollary we have an extension of
Gieseker stability to a t-stability on the bounded derived category $D^b(\Coh X)$. A detailed
discussion of this t-stability is the subject of another article.

%%%%%%%%%%%%%%%%%%%%%%%%%%%%%

\bigskip

The helix generalization (\cite{GK}) of the Beilinson theorem about the basis of the derived
category $D^b(\Coh \PP^n)$ of coherent sheaves on the projective space $\PP^n$ (\cite{BEI})
leads to the remarkable t-stabilities for this category.

Let $\T$ be a linear triangulated category, i.e. for any pair of objects $X,\,Y\in\T$ the direct
sum  $\bigoplus\limits_{j}\Hom^j(X,Y)$ is a   finite dimensional graded vector space. Recall
that an object $E\in\T$ is called exceptional, if $\Hom^\com(E,E)$ is 1-dimension algebra
generated by the identity map. An ordered collection of exceptional objects $(E_0,E_1,\ldots,
E_n)$ is called exceptional, whenever \mbox{$\Hom^\com (E_j,E_i)=0$}\quad $\forall\, i<j$. We
say that the category $\T$ is generated by the exceptional collection, if the smallest full
triangulated subcategory containing all ``linear combinations'' $\bigoplus V_i^\com\otimes E_i$,
where $V_i^\com$ are graded vector spaces, coincides with $\T$. Here for a graded vector space
$V^\com$ and an object $X\in\T$ we denote by $V^\com\otimes X$ the direct sum
$\bigoplus\limits_j V^j\otimes X[-j]$.

For example, the bounded derived category $D^b(\Coh \PP^n)$ is generated by the exceptional
collection $(\O(k),\O(k+1),\ldots, \O(k+n))$, where $k$ is an arbitrary integer number and
$\O(m)$ means the line bundle on $\PP^n$ of degree~$m$.

In this notation we formulate theorem (\cite{GK}) generalizing the Beilinson theorem.

\begin{Thm}\label{BeilHel} Let $\T$ be a linear triangulated category, generated by an
exceptional collection $(E_0,E_1,\ldots, E_n)$. Then for any object $X\in \T$ there exists a
canonical Postnikov system
\[\begin{diagram}[size=1.5em]
&&V_0^\com\otimes E_0&&&&V_1^\com\otimes E_1&&&&&&&&V_n^\com\otimes E_n\\
&\ruTo&&\rdDashto&&\ruTo&&\rdDashto&&&&&&\ruTo&&\rdDashto&&\\
X&&\lTo&&F^1X&&\lTo&&F^2X&\lTo\relax&\cdots&\lTo&F^nX&&\lTo&&0
\end{diagram}
\]
\end{Thm}

Under the assumption of Theorem \ref{BeilHel} let us denote by $\Upsilon_i$ the full
extension-closed subcategory in $\T$, generated by objects $E_i[z]$ with $z\in\ZZ$.

Since the collection $(E_0,E_1,\ldots, E_n)$ is exceptional, the subcategories $\Upsilon_i$
satisfy the property $\Hom^{\le 0}(\Upsilon_j,\Upsilon_i)=0$ whenever $j>i$. Taking into account
Theorem \ref{BeilHel}, we obtain stability data on $\T$. More exactly the following proposition
takes place.

\begin{Prop}\label{Helstab}
Let $\T$ be a linear triangulated category, generated by an exceptional collection
$(E_0,E_1,\ldots, E_n)$. Then $\DST \Delta i \Upsilon$ makes stability data on $\T$, where
$\Delta$ is the natural ordering set $\{0,1,\ldots,n\}$.
\end{Prop}

In the next section we show that \HN\ has the familiar properties  of the Harder--Narasimhan
filtration defined for coherent sheaves in respect to Gieseker or Mumford--Takemoto stability.

\section{Properties of \PGN}\label{PROPERTIES OF PGN}

\noindent  We fix  stability data $\left(\Phi,\{\Pi_\vf\}_{\vf\in\Phi}\right)$ on a triangulated
category $\T$. The important   property of \HN\ is that it is a canonical t-filtration.

\begin{Thm}\label{uniq}
The \HN\ for any object $X\in \T$ is determined up to a unique isomorphism of Postnikov's
systems.
\end{Thm}
\begin{proof}{.} Let
\begin{equation}\label{first HN}
\TFILTR X X F n {\vf} {q} {p}\hspace*{1em}
\end{equation}
and
\begin{equation}\label{second HN}
\TFILTR X X Q m {\psi} {q'\!\!} {p'\!\!}\hspace*{1.5em}
\end{equation}
%\[
%\begin{diagram}[size=1.4em]
%                         &                           &{X}_{{{\psi}}_0}             &         &          &                           &{X}_{{{\psi}}_1}             &         &          &          &      &    &&&{X}_{{{\psi}}_{m}}\\
%                         &\ruTo^{\scriptstyle q'\!\!_0}&                          &\rdDashto&          &\ruTo^{\scriptstyle {q'}_1}&                          &\rdDashto&          &          &      &    &&\ruTo^{\scriptstyle {q'}_{m}}&&\rdDashto\\
%\llap{${X}=$}{Q}^0{X} &                           &\lTo_{\scriptstyle {p'}_1}& &{Q}^1{X}&
%&\lTo_{\scriptstyle {p'}_2}& &{Q}^2{X}&\lTo\relax&\cdots&\lTo&{Q}^{m}{X}&&\lTo_{\scriptstyle
%{p'}_{{m}+1}}&&{Q}^{{m}+1}{X}=0
%\end{diagram}
%\]
 be two \HN s for an object $X\in\T$. We
have show that $n=m$, $\vf_i=\psi_i$ $\forall\,i$, and the identical morphism
$id_X\in\Hom_\T(X,X)$
 induces a unique isomorphism of t-filtrations:
\[
\begin{diagram}[size=2.1em]
   &                        & X_{\vf_0} &         &      &                        &  X_{\vf_1}   &         &     &                &      &                        &  X_{\vf_n}   &         &            \\
   &\ruTo^{\scriptstyle q_0}& \dLine    &\rdDashto&      &\ruTo^{\scriptstyle q_1}&\dLine        &\rdDashto&     &                &      &\ruTo^{\scriptstyle q_n}&\dLine        &\rdDashto&            \\
X  &                        &\lTo       &         & F^1X &                        &\lTo          &         &F^2X &\quad\cdots\quad& F^nX &                        &\lTo          &         &F^{n+1}X        \\
\dTo<{\scriptstyle id_X}>{\scriptstyle \cong}&&\dTo>{\scriptstyle \cong}&&\dTo>{\scriptstyle \cong}%
&&\dTo>{\scriptstyle \cong}&&\dTo>{\scriptstyle \cong}&&\dTo>{\scriptstyle \cong}&&\dTo>{\scriptstyle \cong}&&\dTo>{\scriptstyle \cong}\\
&&X_{\scriptstyle \psi_0}&&&&X_{\scriptstyle \psi_1}&&&&&&X_{\scriptstyle \psi_n}\\
&\ruTo^{\scriptstyle q'_0}&&\rdDashto&&\ruTo^{\scriptstyle q'_1}&&\rdDashto&&&&\ruTo^{\scriptstyle q'_n}&&\rdDashto\\
 X&&\lTo&&Q^1X&&\lTo&&Q^2X&\quad\cdots\quad&Q^nX&&\lTo&&Q^{n+1}X%\llap{\ \ \ \ \ \ $=0$}
\end{diagram}
\]

To do this we prove  some additional properties of  \HN s:
\begin{Prop}\label{hom prop of HN}
Let
\[\TFILT X X F n {\vf}\]
be a \HN\ for $X$. Then
\begin{enumerate}
  \item\label{hom to lesss} $\Hom^{\le 0}(X,\Pi_\vf)=0$ whenever $\vf<\vf_0$,
  \item\label{hom memb to quot} $\Hom^{\le 0}(F^iX,\Pi_\vf)=0$ whenever $\vf\le \vf_i$,
  \item\label{hom from great} $\Hom^{\le0}(\Pi_\psi,X)=0$ whenever $\psi>\vf_n$;
  \item\label{whole hom} if
  \[\TFILT Y Y F m {\psi}\]
  is a \HN\ for $Y$ such that $\vf_n<\psi_0$, then $\Hom^{\le0}(Y,X)=0$.
\end{enumerate}
\end{Prop}

\begin{proof}{ of the proposition} is easy obtained via the application of functors
$\Hom(\cdot,\Pi_\vf)$,\linebreak $\Hom(\Pi_\psi,\cdot)$, and $\Hom(Y,\cdot)$ to each of the
triangle
\[\TR{F^{i+1}X}{F^iX}{X_{\vf_i}}\qquad.\]
One should  notice that $F^nX\in \Pi_{\vf_n}$. We leave the details   to the reader.\end{proof}

Returning to the proof of the theorem, consider the first triangle of \HN\ \eqref{first HN} for
$X$ and any $Y\in \Pi_{\vf_0}$. Applying the functor $\Hom(\cdot,Y)$ to the triangle, we have
\[ \Hom^{-1}(F^1X,Y)\To\Hom^0(X_{\vf_0},Y)\xTo{h_Y}\Hom^0(X,Y)\To\Hom^0(F^1X,Y) \ .\]
 It follows from Proposition \ref{hom prop of HN} (\ref{hom memb to quot}), that the very left and
  right terms  of the
above  exact sequence are zero. Hence,
\begin{equation}\label{repr.obj}
  \Hom^0(X_{\vf_0},Y)\cong \Hom^0(X,Y)\quad\forall\,Y\in\Pi_{\vf_0}.
\end{equation}
In other words, the object $X_{\vf_0}$ represents the functor $\Hom^0(X,\cdot):\Pi_{\vf_0}\To
\Vect$.

Let us substitute $X_{\vf_0}$ for $Y$ in \eqref{repr.obj} and denote the morphism
$h_{X_{\vf_0}}(id_{X_{\vf_0}})\in\Hom^0(X,X_{\vf_0})$ by $q_0$. It can be shown in the usual way
(see, for example, \cite[Lemma IV.4.5]{GM}) that the \tr\
\[
\TRM{F^1X}X{X_{\vf_0}}{}{q_0}\quad\ , \] where $X_{\vf_0}\in\Pi_{\vf_0}$ and
$\Hom^{\le0}(F^1X,X_{\vf_0})=0$, is determined up to a unique isomorphism of triangles.

Note that $q_0\ne0$. Therefore, $\Hom^0(X,\Pi_{\vf_0})\ne0$. Further, starting with the first
triangle of \HN\ \eqref{second HN}, we get $\Hom^0(X,X_{\psi_0})\ne0$. Applying now Proposition
\ref{hom prop of HN} (\ref{hom to lesss}) to $X$, $\Pi_{\vf_0}$ and $\Pi_{\psi_0}$, we conclude
that $\vf_0=\psi_0$. Moreover, it follows from the uniqueness  of the   object representing a
functor, that the identical map $id_{X}\in \Hom_\T(X,X)$  extends  to a canonical isomorphism of
triangles:
\[
\begin{diagram}[size=1.5em]
X_{\vf_0}&           & &\rTo^\sim_{\scriptstyle \exists!}& &            &X_{\psi_0}\\
\dDashto &\luTo_{\scriptstyle q_0}& &                    & &\ruTo_{\scriptstyle q'_0}&\dDashto  \\
         &           &X&\rTo^\sim_{\scriptstyle id_X}    &X&            &          \\
         &\ruTo      & &                    & &\luTo       &          \\
F^1X     &           & &\rTo^\sim_{\scriptstyle \exists!}& &            &Q^1X
\end{diagram}
\]
The rest of the proof is done  by induction.\end{proof}

The last proposition that we prove in this section  deals with properties of a general\\
t-filtration.

\begin{Prop}\label{atom. and glu. t-filtr}
\begin{enumerate}
\item Let $X$ has a t-filtration $X\rightsquigarrow (Y_0,Y_1,\ldots,Y_n)$ and  each $Y_s$ has one
 too\\
     $Y_s\rightsquigarrow(X_{ (s,0)},\ldots,X_{ (s,k_s)})$ $\forall\,s$. Then one can construct
     a combined  t-filtration
\[X\rightsquigarrow(X_{(0,0)},\ldots,X_{ (0,k_0)},X_{ (1,0)},\ldots,X_{(1,k_1)},\ldots,
    X_{(n,0)},\ldots,X_{(n,k_n)}).\]
    \item Let
\[X\rightsquigarrow(X_{(0,0)},\ldots,X_{ (0,k_0)},X_{ (1,0)},\ldots,X_{ (1,k_1)},\ldots,
    X_{ (n,0)},\ldots,X_{ (n,k_n)})\] be a t-filtration. Then there exist  t-filtrations
    $X\rightsquigarrow (Y_0,Y_1,\ldots,Y_n)$ and\\
     $Y_s\rightsquigarrow(X_{ (s,0)},\ldots,X_{ (s,k_s)})$    for each $s$.

\item Let $X\rightsquigarrow(X_0,X_1,\ldots,X_n)$ and   $Y\rightsquigarrow(Y_0,Y_1,\ldots,Y_m)$ be
t-filtrations. Then for each shuffle permutation of the quotients, i.e.  linear ordering of the
set
\[\{Z_s\}_{s=0,\ldots,n+m+1}=\{X_i\}_{i=0,1\ldots,n}\bigsqcup \{Y_j\}_{j=0,\ldots,m}\]
that respects the initial  orders on  $\{X_i\}_{i=0,1\ldots,n}$ and $\{Y_j\}_{j=0,\ldots,m}$
there exists a t-filtration $X\oplus Y\rightsquigarrow(Z_0,Z_1,\ldots,Z_{n+m+1})$.
\end{enumerate}
\end{Prop}
\begin{proof}{.} We start with the second statement.
    Let t-filtration for $X$ in the   statement be determined by the sequence of \tr s:
\begin{align*}
 &\TRM{F^{(i,j+1)}X}{F^{(i,j)}X}{X_{(i,j)}}{p_{(i,j+1)}}{}\qquad,\qquad\text{where } i=0,\ldots,n;\ j=0,\ldots,k_i-1;\text{
 and}\\[2ex]
 &\TRM{F^{(i+1,0)}X}{F^{(i,k_i)}X}{X_{(i,k_i)}}{p_{(i+1,0)}}{}\qquad,\qquad\text{where }
 i=0,\ldots,n;
\end{align*}
$ F^{(0,0)}X=X$;  and $F^{(n+1,0)}X=0$.

Denote the composition $p_{(i+1,0)}p_{(i,k_i)}\cdots p_{(i,1)}: F^{(i+1,0)}X\To F^{(i,0)}X$ by
$p_i$. In these notations we have got  a new filtration for $X$ that consists of   triangles:
\[
\TRM{F^{(i+1,0)}X}{F^{(i,0)}X}{Y_i}{p_i}{}\qquad.
\]
Therefore, we need to show that for any $i=0,\ldots,n$ there exists a t-filtration
$$Y_i\rightsquigarrow(X_{(i,0)},\ldots, X_{(i,k_i)}).$$ It is clear that there is a t-filtration
\[
    F^{(i,0)}X\rightsquigarrow (X_{{(i,0)}},\ldots,X_{ {(i,k_i)}},X_{ {(i+1,0)}},
    \ldots,X_{ {(n,k_n)}})
\]
and the statement    follows from the lemma:
\begin{Lem}\label{glue of t-filtr}
    Let
\[
    \begin{diagram}[size=1.4em]
      &          & Z_0      &         &    &          & Z_1      &         &    &          &      &     &    &     &Z_m           &         &\\
      &\ruTo     &          &\rdDashto&    &\ruTo     &          &\rdDashto&    &          &      &     &    &\ruTo&              &\rdDashto&\\
Z=F^0Z&          &\lTo_{\scriptstyle h_1}&         &F^1Z&          &\lTo_{\scriptstyle h_2}&
&F^2Z&\lTo\relax&\cdots&\lTo &F^mZ&     &\lTo_{\scriptstyle h_{m+1}}&         &F^{m+1}Z=0
    \end{diagram}
\]
be a t-filtration for an object $Z$. For each $s$, $1\le s\le m+1$ denote by $c_s$ the
composition $c_s=h_sh_{s-1}\cdots h_1: F^sZ\To Z$ and consider   \tr s
\[\TRM{F^sZ}Z{Q_s}{c_s}{}\quad.\]
Then for each $s$ there exists a t-filtration $Q_s\rightsquigarrow(Z_0,Z_1,\ldots, Z_{s-1})$.
\end{Lem}
\begin{proof}{ of the lemma} is by induction on s. For $s=1$ there is nothing to prove. In the
general case consider the diagram of \tr s:
\[
\begin{diagram}[size=1.4em]
     &    &Z_s         &          & \\
     &    &\uTo        &          &  \\
Q_s[-1]&\rTo&F^sZ        &\rTo^{\scriptstyle c_s}&Z\,. \\
     &    &\uTo<{\scriptstyle h_{s+1}}&          &   \\
     &    &F^{s+1}Z    &          &
\end{diagram}
\]
Applying  the octahedron axiom to the diagram
\[
\begin{diagram}[size=1.4em]
     &    &Z_s         &\rTo          &Q_{s+1} \\
     &\ruTo    &\uTo        &          &\uTo  \\
Q_s[-1]&\rTo&F^sZ        &\rTo^{\scriptstyle c_s}&Z\,, \\
     &    &\uTo<{\scriptstyle h_{s+1}}&\ruTo_{\scriptstyle c_{s+1}}          &   \\
 &    &F^{s+1}Z    &          &\end{diagram}\]
  we obtain \tr s:
\[
\TRM{F^{s+1}Z}{Z}{Q_{s+1}}{c_{s+1}}{}\qquad,\qquad \TR{Z_s}{Q_{s+1}}{Q_s}\ \ .
\]
By the induction hypothesis, $Q_s\rightsquigarrow(Z_0,\ldots, Z_{s-1})$. In particular, we have
the first \tr\ of a t-filtration for $Q_s$
\[
\TR{F^1Q_s}{Q_s}{Z_0}\] and $F^1Q_s\rightsquigarrow (Z_1,\ldots,Z_{s-1})$ the rest of the
t-filtration. We shall construct the filtration
$Q_{s+1}\rightsquigarrow(Z_0,Z_1,\ldots,Z_{s-1},Z_s)$.

First we apply  to the diagram  of \tr s
\[
\begin{diagram}[size=1.4em]
     &    &Z_0         &           &  \\
     &     &\uTo        &          &  \\
Q_{s+1}&\rTo&Q_s        &\rTo &Z_s[1] \\
     &    &\uTo &           &   \\
 &    &F^{1}Q_s    &          &\end{diagram}\]
 the octahedron axiom:
\[
\begin{diagram}[size=1.4em]
     &    &Z_0         &\rTo          &F^1Q_{s+1}[1] \\
     &\ruTo    &\uTo        &          &\uTo  \\
Q_{s+1}&\rTo&Q_s        &\rTo&Z_s[1]\,, \\
     &    &\uTo &\ruTo           &   \\
 &    &F^{ 1}Q_s    &          &\end{diagram}\]
and get     $F^1Q_{s+1}$ and \tr s
\[
\TR{Z_s}{F^1Q_{s+1}}{F^1Q_s}\ \ ,\qquad \TR{F^1Q_{s+1}}{Q_{s+1}}{Z_0}\ \ \ \ \ .
\]
  By the induction hypothesis, it follows from the first triangle that  there exists a
  t-filtration
$F^1Q_{s+1}\rightsquigarrow(Z_1,\ldots,Z_{s-1},Z_s)$. Therefore, one can  include the second
triangle into the t-filtration
\[
    \begin{diagram}[size=1.4em]
      &          & Z_0      &         &    &          & Z_1      &         &    &          &      &     &    &     &Z_s           &         &\\
      &\ruTo     &          &\rdDashto&    &\ruTo     &          &\rdDashto&    &          &      &     &    &\ruTo&              &\rdDashto&\\
Q_{s+1}=F^0Q_{s+1}&          &\lTo &         &F^1Q_{s+1}&          &\lTo & &F^2Q_{s+1}&
\lTo\relax&\cdots&\lTo &F^sQ_{s+1}& &\lTo &         &F^{s+1}Q{s+1}=0\,.
    \end{diagram}
\]
This completes the proof of the lemma and the second statement of the
proposition.\end{proof}

Arguing very similarly, one can prove the first statement of the proposition. To prove the last
one, consider the t-filtrations:
\begin{align*}
 &\TF X X F n  {q} \qquad\ \ ,\\[1ex]
 &\TF Y Y F m  {p} \qquad\ \ .
\end{align*}
Suppose,
\[
    (Z_0,\ldots,Z_{n+m+1})=(X_0,\ldots,X_{i_1},Y_0,\ldots,Y_{j_1},X_{i_1+1},\ldots, ).\]
One constructs the t-filtration
\[
\begin{diagram}[size=1.4em]
&&Z_0&&&&Z_{1}&&&&Z_{2}\\
&\ruTo^{\scriptstyle h_0}&&\rdDashto&&\ruTo^{\scriptstyle h_1}&&\rdDashto&&\ruTo^{\scriptstyle h_2}&&\rdDashto\\
X\oplus Y &&\lTo& &F^1(X\oplus Y)& &\lTo& &F^2(X\oplus Y)&&\lTo&&F^3(X\oplus
Y)&\lTo\relax&\cdots
\end{diagram}
\]
in the following way:
\begin{align*}
  &h_0=q_0\oplus 0, \,\ldots,\, h_{i_1}=q_{i_1}\oplus 0;\\
  &F^1(X\oplus Y)=F^1X\oplus Y,\,\ldots ,\, F^{i_1}(X\oplus Y)=F^{i_1}X\oplus Y;\\
  &h_{i_1+1}=0\oplus p_0, \,\ldots,\, h_{i_1+j_1}=0\oplus p_{j_1};\\
  &F^{i_1+1}(X\oplus Y)=F^{i_1}X\oplus F^1Y,\,\ldots ,\,
   F^{i_1+j_1}(X\oplus Y)=F^{i_1}X\oplus F^{j_1}Y;
\end{align*}
and so on.\end{proof}

\begin{Cor}\label{direct summands} Let $X\in \T$ be a semistable object of slope $\vf$. If
$X=X^1\oplus X^2\oplus\cdots\oplus X^k$, then each summand $X^i$ is semistable of slope $\vf$ as
well.
\end{Cor}
\begin{proof}{.} It is sufficient to prove the statement for $k=2$. Consider the \HN s for the
direct summands
\[X^1\rightsquigarrow (X^1_{\vf_0},X^1_{\vf_1},\ldots, X^1_{\vf_n}),\quad
  X^2\rightsquigarrow (X^2_{\psi_0},X^2_{\psi_1},\ldots, X^2_{\psi_m}).\]
Ordering the quotients in compliance with the order on $\Phi$, unless $n=m=0$ and $\vf_0=\psi_0$
one can construct the nontrivial \HN\ for $X$ (see Proposition \ref{atom. and glu. t-filtr}) .

But $X$ has only the trivial \HN\ because it is semistable (Theorem \ref{uniq}). This concludes
the proof.
\end{proof}

\section{Connections between t-stabilities and t-structures}\label{ORDER}

\noindent  It was shown in Section \ref{DEFINITION} that a bounded t-structure on a triangulated
category induces stability data (Lemma \ref{stability by t-str}). On the other hand a
t-stability gives a family or collection of t-structures. Before formulation  more exact
statement, let us introduce some convenient notations.

 For a subset $S$ of a triangulated category $\T$ we denote by $\LA
S\RA$ the minimal full extension-closed subcategory of $\T$, containing the subset $S$. Note
that if we have a t-stability $\DST \Phi\vf \Pi$ on $\T$ and $\Psi\subset\Phi$, then $\LA
\Pi_\vf|\ \vf\in\Psi\RA$ consists of objects whose have \HN\
$X\rightsquigarrow(X_{\psi_0},\ldots,X_{\psi_k})$ with $\psi_i\in\Psi$.

\begin{Lem}\label{t-structures by t-stability}
    Let $\DST \Phi \vf \Pi$ makes stability data on a triangulated category $\T$. Suppose there
 exists     a decomposition of the set
    $\Phi$ in two disjoint parts: $\Phi=\Phi_- \sqcup \Phi_+$ such that for any
    $\vf_-\in\Phi_-$ and $\vf_+\in\Phi_+$ we have $\vf_-<\vf_+$. Then the pair of subcategories
\[ \T^{\ge 0}=\LA\Pi_\vf|\ \vf\in \tau(\Phi_-)\RA,\quad \T^{\le 0}=\LA\Pi_\vf|\ \vf\in \Phi_+ \RA
\]  defines a t-structure on $\T$ (recall that $\tau$ is an automorphism of the linearly ordered set
$\Phi$, that corresponds to the shift of $\T$).
\end{Lem}

\begin{Cor}\label{t-struct by vf}
Each element $\vf\in\Phi$ gives a
 t-structure
 \[ \T^{\ge 0}_\vf=\LA\Pi_\psi|\ \psi\le \tau(\vf)\RA,\quad \T^{\le 0}_\vf=\LA\Pi_\psi|\
 \psi>\vf \RA.
\]\end{Cor}
\begin{proof}{ of the lemma.}
 Firstly let us recall that by the definition of t-stability $\Pi_\vf[-1]=\Pi_{\tau^{-1}(\vf)}$ and
 $\tau(\vf)\ge \vf$, whence the first axiom of a t-structure is valid:
\[\T^{\le 0}\subset \T^{\le 0}[-1]\quad\text{ and }\quad \T^{\ge 0}\supset \T^{\ge 0}[-1].\]
Besides
 $\Hom^{\le0}(\Pi_\vf, \Pi_\psi)=0$ whenever $\vf>\psi$. Therefore the axiom
 \[\Hom^{\le0}(\T^{\le 0},\T^{\ge 0}[-1])=0\] of a t-structure follows from the statement:
 \begin{quote}\itshape
    Let  $S,\,S'$ be   subsets of objects of $\T$ such that $\Hom^{\le0}(X,Y)=0$ for each
    $X\in S$ and $Y\in S'$. Then $\Hom^{\le0}(\LA S\RA,\LA S'\RA)=0$,
 \end{quote}
 which is evidently true.

To verify the last axiom consider \HN\ for an object $0\ne X\in\T$:
\[
\TFILTR X X F n {\vf}{q}{p}\qquad\quad\ .
\]
If $\vf_n\in \Phi_-$ or $\vf_0\in \Phi_+$, then the \tr\ of the axiom is trivial. Suppose  there
exists $k$ such that $\vf_k\in \Phi_-$ and $\vf_{k+1}\in \Phi_+$. Completing the morphism
$X_{\le0}=F^{k+1}X\xTo{p_{k+1}p_k\cdots p_1}X$ by a cone $X_{\ge1}$, we get the needed \tr
\[\TR {X_{\le0}}X{X_{\ge1}}\]
of the axiom (see Proposition \ref{atom. and glu. t-filtr}).
\end{proof}

\medskip

\noindent {\bf Remark.} {\itshape Note that the t-structure, constructed above, is bounded if
and only if }
\[
\Phi=\left(\bigcup_{n\in\ZZ_{\ge0}}\tau^{-n}(\Phi_+)\right)\quad\text{and}\quad
\Phi=\left(\bigcup_{n\in\ZZ_{\ge0}}\tau^{n}(\Phi_-)\right).
\]

As a result, we have got the connection between t-structures and t-stabilities on a given
triangulated category.

It is easy to observe, that a given t-structure can be induced by several stability data.

Let us consider a t-stability $\DST \Phi\vf\Pi$ on $\T$ such that $\Phi$ is a disjoin union
$\Phi=\coprod\limits_{\psi\in\Psi}\Phi_\psi$ of nonempty subsets  with the following properties:
\begin{enumerate}
\item if  $\vf_i\in\Phi_{\psi_i}$ ($i=1,2$) and  $\vf_1<\vf_2$, then for any
      $\vf'_i\in \Phi_{\psi_i}$ the inequality $\vf'_1<\vf'_2$ holds.
\item $\forall\, \psi\in \Psi$ \ $\exists\, \psi'\in\Psi$ such that
$\tau(\Phi_\psi)=\Phi_{\psi'}$.
\end{enumerate}

The properties of the disjoin union allow to determine an order on $\Psi$ and $\bar\tau\in
\Aut(\Psi)$, corresponding to $\tau \in\Aut (\Phi)$ in the obvious way.

Further we define $P_\psi$ for $\psi\in\Psi$ as $P_\psi=\LA \Pi_\vf|\ \vf\in \Phi_\psi\RA$. The
fact that $\DST \Psi \psi P$ is a t-stability follows immediately from the definition and
Proposition \ref{atom. and glu. t-filtr}.

We have constructed another t-stability $\DST \Psi\psi P$, but a t-structure induced by a
decomposition of $\Psi$, can be obtained also from a decomposition of $\Phi$. The t-stability
$\Phi$ induces all t-structures that $\Psi$ does, and more. We would like to say  that the
stability data $\Phi$ are finer then $\Psi$, and $\Psi$ are weaker. Generalizing this
construction, let us formulate the definition.

\begin{Defn}\label{fine and weak}
    Let $\DST \Phi \vf \Pi$ and $\DST \Psi \psi P$ be t-stabilities on a triangulated category
    $\T$ with automorphisms $\tau_\Phi$ and $\tau_\Psi$, corresponding to the shift on $T$. We
    say that $\Phi$ is finer then $\Psi$ (and $\Psi$ is weaker then $\Phi$) and denote this by
    $\Phi\preceq\Psi$ if there exists a surjection $r:\Phi\To\Psi$ such that
\begin{enumerate}
  \item $r\tau_\Phi=\tau_\Psi r$;
  \item $\vf'>\vf''$ $\Leftrightarrow$ $r(\vf')\ge r(\vf'')$;
  \item $\forall\,\psi\in\Psi$\quad $P_\psi=\LA \Pi_\vf|\ \vf\in r^{-1}(\psi)\RA$.
\end{enumerate}
\end{Defn}

Clearly this gives a partial order on the set of all t-stabilities on a given triangulated
category. Therefore, some of stability data can be minimal w.r.t the order.  They seem to
contain  the maximal information about t-structures. We call such t-stabilities the finest.

Now we give conditions of comparability and being the  finest  for stability data.

\begin{Prop}\label{comp. and fine. cond}
     Let $\DSTm \Phi \vf \Pi$ and $\DSTm \Psi \psi P$ make stability data on $\T$. Then

1. $\Phi\preceq\Psi$ if and only if
\begin{enumerate}
  \item[(i)] any $\Phi$-semistable object is $\Psi$-semistable,
  \item[(ii)] for $\Pi_{\vf_i}\subset P_{\psi_i}$ ($i=1,2$) the condition $\vf_1<\vf_2$ implies
  $\psi_1\le \psi_2$,
  \item[(iii)] if $\Pi_\vf\subset P_\psi$, then $\Pi_{\tau_\Phi(\vf)}\subset
  P_{\tau_\Psi(\psi)}$.
\end{enumerate}

2. Suppose that for each $\vf\in\Phi$ we have
\[ \forall\,X,Y\in\Pi_\vf\qquad\Hom_\T(X,Y)\ne 0 \ \text{ and } \ \Hom_\T(Y,X)\ne 0.\]
Then the t-stability $\DSTm \Phi\vf \Pi$ is the finest.
\end{Prop}

\begin{proof}{.}
The fact that the assumption $\Phi \preceq\Psi$ implies conditions (i)--(iii) follows
immediately  from the definition. To prove the converse  it is sufficient to show that each
semistable subcategory $P_\psi$ is generated by some collection of subcategories $\Pi_\vf$. In
other words, if $X_\psi\in P_\psi$ has a \HN\
$X_\psi\rightsquigarrow(X_{\vf_0},\ldots,X_{\vf_n})$ w.r.t t-stability $\Phi$, then
$\Pi_{\vf_i}\subset P_\psi$ for each $i$.

By condition (i) any $X_{\vf_i}$ is $\Psi$-semistable. Therefore, we can find $\psi_i\in\Psi$
such that $X_{\vf_i}\in P_{\psi_i}$. Since $\Hom^0(X_\psi, X_{\vf_0})\ne 0$ and $X_{\vf_0}\in
P_{\psi_0}$, we have $\psi\le \psi_0$ (Proposition \ref{hom prop of HN}). On the other hand,
$\Hom^0(X_{\vf_n},X_\psi)\ne 0$, consequently, $\psi_n\le\psi$. Besides,
$\vf_0<\vf_1<\cdots<\vf_n$, so (by condition (ii)) $\psi_0\le\psi_1\le\cdots\le \psi_n$. This is
possible only if $\psi_0=\psi_n$.

To verify the second statement of the proposition, we note that if $\DST \Phi \vf \Pi$ is not
the finest, then there exists  some $\vf\in \Phi$ and nonempty  subcategories $\Pi_-$, $\Pi_+$
such that $\Pi_\vf=\LA \Pi_-,\Pi_+\RA$, where    $\Hom^{\le 0}(\Pi_+,\Pi_-)=0$. This contradicts
the assumption of the statement.\end{proof}

{\bf Remark.} {\itshape As we shall see later, a pair of t-stabilities can satisfy the condition
(i) and be noncomparable. The second condition is significant\/}.

\medskip

As the conclusion of the section we formulate a   statement that we shall  use later for the
classification of bounded t-structures.

\begin{Prop}\label{class t-struct}
    Let $\T$ be a triangulated category with the property that for any t-stability $\Psi$ there
    exists the finest t-stability $\Phi\preceq\Psi$. Then for any bounded t-structure
    $(\T^{\le0},\T^{\ge0})$ there exists a t-stability $\Phi=\DST \Phi \vf \Pi$ and a
    decomposition $\Phi=\Phi_-\sqcup\Phi_+$ in two disjoint parts such that:
\[\T^{\le 0}=\LA \Pi_\vf|\ \vf\in \Phi_+\RA,\qquad
\T^{\ge 0}=\LA \Pi_\vf|\ \vf\in \tau(\Phi_-)\RA.\] Without loss of generality we can assume that
$\Phi$ is one of the finest t-stabilities on $\T$.
\end{Prop}
\begin{proof}{} follows directly from Lemmas \ref{stability by t-str} and
\ref{t-structures by t-stability}, and assumptions of the proposition.\end{proof}

\section{Stability data for ${\boldsymbol {D}^{\boldsymbol b}\boldsymbol {(}
{\boldsymbol \Coh}{\boldsymbol \PP}^{\boldsymbol 1}{\boldsymbol )}}$}\label{PROJECTIVE LINE}

\noindent
Here we study t-stabilities and bounded t-structures on the bounded
derived
category $D^b(\Coh \PP^1)$ of coherent sheaves on the projective line. At first  we describe two
types of the finest t-stabilities for $D^b(\Coh\PP^1)$: standard and exceptional. Then we show
that any t-stability for $D^b(\Coh\PP^1)$ is weaker than either standard one or exceptional one.
This
enables us to classify all% At last using the finest stability data, we give the whole list of
bounded t-structures for $D^b(\Coh\PP^1)$.

In this   section  we identify sheaves with 0-complexes of $D^b(\Coh \PP^1)$.

\subsection{Standard t-stability}\label{Standard t-stability}
It follows from Lemma \ref{stability by t-str}, that $\left(\ZZ,
\{(\Coh\PP^1)[i]\}_{i\in\ZZ}\right)$ makes stability data for  $D^b(\Coh \PP^1)$. Obviously,
they are not the finest. To construct the finest t-stability $\DST \M \mu \Pi \preceq \left(\ZZ,
\{(\Coh\PP^1)[i]\}_{i\in\ZZ}\right)$, recall that each sheaf $F$ on $\PP^1$ decomposes in a
finite direct sum
\[
    F=\Bigl(\bigoplus_j \Xi_{x_j}\Bigr)\oplus \Bigl(\bigoplus_i \O(n_i)\Bigr),
\]
where $\Xi_{x_j}$ is a torsion sheaf concentrated at the point $x_j$, and $\O(n_i)$ is an
invertible sheaf of degree $n_i$. Let us define semistable subcategories as
\begin{equation}\label{standard ssscat}
\Pi_{(i,x)}=\LA \O_x[i]\RA\quad\text{and}\quad\Pi_{(j,n)}=\LA \O(n)[j]\RA
\end{equation}
with $x\in\PP^1$, and $n,i,j\in\ZZ$. This yields the set of slopes $\M=\ZZ\x (\ZZ\sqcup\PP^1)$.

Now we have to introduce a linear order on $\M$ in a way
\[ \mu>\nu\ \Rightarrow\ \Hom^{\le 0}(\Pi_\mu,\Pi_\nu)=0.\]
This condition implies:
\begin{align*}
  & (i,\alpha)>(j,\beta)&&\text{for } i>j,\ \alpha,\,\beta\in \ZZ\sqcup \PP^1;\\
  & (i,n)>(i,m)&&\text{for } n>m\in\ZZ,\ i\in\ZZ;\\
  & (i,x)>(i,n)&&\text{for } x\in\PP^1,\ i,\,n\in \ZZ.
\end{align*}
Besides, for each $i,\,j,\,q\in\ZZ$\quad $\Hom^q(\O_x[i],\O_y[j])=0$ unless $x=y$. Thus,
defining a linear order on $\PP^1$ in arbitrary way, we get the linearly ordered set $\M$.

To prove that $\DST \M \mu \Pi$ is a t-stability it is sufficient to verify that each nonzero
object has \HN. This   directly follows from Proposition \ref{atom. and glu. t-filtr} and the
fact that \[\Coh\PP^1[i]=\LA\Pi_{(i,\alpha)}|\ \alpha\in \ZZ\sqcup \PP^1\RA.\]

Finally, note that all the  semistable subcategories $\Pi_{(i,\alpha)}$ ($\alpha\in \ZZ\sqcup
\PP^1$) satisfy the assumption of Proposition \ref{comp. and fine. cond}\,(2). Therefore  the
stability data $\DST \M \mu \Pi$ is the finest one.

Notice, that the t-stability $\DST \M \mu \Pi$ depends on an order on the slope  set $\M$. Hence
we have got various finest   incomparable   t-stabilities with the same semistable subcategories
   and slopes. Each of them we shall call   standard.

\subsection{Exceptional t-stability}\label{Exceptional t-stability}
The construction of the next type of stability data is based on Proposition \ref{Helstab}. The
category $D^b(\Coh\PP^1)$ is generated by the exceptional collection \mbox{$(\O(k), \O(k+1))$}
($k$ is an arbitrary fixed integer number). Therefore we have stability data $\left(\{0,1\},
\{\Pi^k_i\}_{i\in\{0,1\}}\right)$, for each  integer number $k$, where
\begin{align*}
  &\Pi^k_0=\LA \O(k)[j]|\ j\in\ZZ\RA,\\
  &\Pi^k_1=\LA \O(k+1)[i]|\ i\in\ZZ\RA.
\end{align*}

\medskip

Since $\Hom^{\le0}(\O(n)[i],\O(n)[j])=0$ for  $n,\,i,\,j\in\ZZ$ whenever $i>j$, the t-stability
defined above is not finest. To refine it we consider semistable subcategories and   slopes
\[ \Pi^k_{(i,0)}=\LA \O(k)[i]\RA,\quad \Pi^k_{(i,1)}=\LA \O(k+1)[i]\RA,\quad
(i,0);\,(i,1)\in\E=\ZZ\x\{0,1\}.\]
We are going to refine the stability data $\left(\{0,1\},
\{\Pi^k_i\}_{i\in\{0,1\}}\right)$ to some finest
$\left(\E,\set{\Pi^k_\eps}_{\eps\in\E}\right)$.
%$\left(\E,\set{\Pi^k_\eps}_{\eps\in\E}\right)\preceq $\left(\{0,1\},\{\Pi^k_i\}_{i\in\{0,1\}}\right)$.
A~linear order on the set $\E$ and automorphism $\tau\in\Aut\E$   should be coordinated with the
order and the automorphism above. Therefore, we come to
\begin{align*}
  &(i,0)<(j,1);\\
  &(i,0)<(i+1,0),\ (i,1)<(i+1,1);\\
  &\tau(i,0)=(i+1,0),\  \tau(i,1)=(i+1,1)
\end{align*} for arbitrary $i,j\in\ZZ$.

Finally we need to obtain a finite \HN\ w.r.t. $\E$ for each nonzero object $X$ of the derived
category.  Note that homological dimension of $\Coh\PP^1$ is 1. We shall  use the following
result well-known among specialists, but it seems there is no reference.

\begin{Prop}\label{dh=1}
    If homological dimension of an abelian category $\A$ equals 1, then each object $X$ of $D^b(\A)$
    is isomorphic to a finite direct sum $X=\bigoplus\limits_i A_i[-i]$, where $A_i\in\A$ and
    $A_i[-i]$ denotes   $i$-complex, i.e., a complex with the unique nonzero term  $A_i$
    at the place  $i$.
\end{Prop}

 For convenience of the reader we provide the proof.

\begin{proof}{.} Suppose $X\in D^b(\A)$ is represented by a complex
\[C^{\com}=\quad C^n\xTo{d_n}C^{n+1}\xTo{d_{n+1}}\cdots\To C^{m-1}\xTo{d_{m-1}} C^m\]
such that $H^i(C^\com)=0$ for $i<s$ and $H^s(C^\com)\ne0$. The exact sequence of complexes
\[
\begin{diagram}[size=1.5em]
K^\com&=\qquad&C^n&\rTo^{\scriptstyle d_n}&C^{n+1}&\rTo\relax&\cdots&\rTo&C^{s-1}&\rTo^{\scriptstyle d_{s-1}}&\Im d_{s-1}&\rTo&0\\
\dInto&&\dTo>{\scriptstyle id}   &&\dTo>{\scriptstyle id}  &&&&\dTo>{\scriptstyle id}   &&     \dInto &&\dTo  \\
C^\com&=\qquad&C^n&\rTo^{\scriptstyle d_n}&C^{n+1}&\rTo\relax&\cdots&\rTo&C^{s-1}&\rTo^{\scriptstyle d_{s-1}}&C^s&\rTo^{\scriptstyle d_s}&C^{s+1}&\rTo\relax&\cdots\\
\dOnto&& &&&&&&\dTo &&\dOnto &&\dTo>{\scriptstyle id}       \\
B^\com&=\qquad&&&&&&&0&\rTo&B^s&\rTo^{\scriptstyle \bar d_s}&C^{s+1}&\rTo\relax&\cdots
\end{diagram}
\]
determines the \tr\ in $D^b(\A)$
\[\TR YXZ\]
with  $Y\simeq K^\com$ and $Z\simeq B^\com$. Since the complex $K^\com$ is acyclic, $Y\simeq 0$
(see \cite{GM}), i.e. $X\simeq Z\simeq B^\com$.

Now we shall prove the proposition by induction on the number $h(B^\com)$ of nonzero
cohomologies of the complex $B^\com$.

The base of induction ($h=1$) was proved in \cite{GM}. Further, let us complete the inclusion
$(\Ker\bar d_s)[-s]\hookrightarrow B^\com$ to a \tr
\begin{equation}\label{splitting triangle}
  \TR{(\Ker\bar d_s)[-s]}{B^\com}{D^\com}\qquad\quad.
\end{equation}
Since $h(D^\com)<h(B^\com)$, we can use the induction hypothesis: $D^\com\simeq
\bigoplus\limits_{i>s}A_i[-i]$ with $A_i\in\A$. Applying to the triangle the functor
$\Hom(\cdot,(\Ker\bar d_s)[-s])$, we get exact sequence
\[\begin{split}
& \Hom^0(B^\com,(\Ker\bar d_s)[-s])\To\Hom^0((\Ker\bar d_s)[-s],(\Ker\bar d_s)[-s])\To
\\\To&
\bigoplus_{i>s}\Hom^1( A_i[-i],(\Ker\bar d_s)[-s]). \end{split}
\] It follows from the
inequality $i>s$   that $q=1+i-s\ge2$ and we have
\[\Hom^1( A_i[-i],(\Ker\bar d_s)[-s])=\Ext^q_\A(A_i,\Ker\bar d_s).\]
 On the other hand, homological dimension of $\A$ equals 1, so
$\bigoplus\limits_{i>s}\Hom^1( A_i[-i],(\Ker\bar d_s)[-s])=0$ and the triangle \eqref{splitting
triangle} splits. This completes the proof.\end{proof}

In addition, as we already marked, a sheaf on $\PP^1$ is direct sum of invertible sheaves and
torsion sheaves with support at   point. Besides,   the shift of the derived category acts on
the set of \tr s. Therefore, according to Proposition \ref{atom. and glu. t-filtr},  it is
sufficient to construct \HN\ for an invertible sheaf $\O(n)$ and a torsion sheaf $\Xi_x$
concentrated at a point $x$. The needed  HN-filtrations  are obtained from the exact
sequences\footnote{When $X$ is an object of an additive category we use notation $mX$ for the
object $X^{\oplus m}$.}
\begin{align*}
 &\seq {(n-k-1)\O(k)}{(n-k)\O(k+1)}{\O(n)},&&\text{if }n>k+1,\\
 &\seq{\O(n)}{(k-n+1)\O(k)}{(k-n)\O(k+1)},&&\text{if }n<k,\\
 &\seq{d\O(k)}{d\O(k+1)}{\Xi_x}&&\text{where }d=\deg \Xi_x ,
\end{align*}
and they have  the form
\begin{align}
 &\TR{(n-k)\O(k+1)}{\O(n)}{(n-k-1)\O(k)[1]}&\text{if }n>k+1,\label{triangles for P1-1}
\end{align} \begin{align}&\TR{(k-n)\O(k+1)[-1]}{\O(n)}{(k-n+1)\O(k)}&\text{if }n<k,\label{triangles for P1-2}\\
 &\TR{d\O(k+1)}{\Xi_x}{d\O(k)[1]}\qquad\,.&\label{triangles for P1-3}
\end{align}

Thus, we have got the finest stability data $\left(\E,\set{\Pi^k_\eps}_{\eps\in\E}\right)$ (see
Prop. \ref{comp. and fine. cond}).

As we saw before, changing an order of the slopes set, one can get another t-stability. Let us
research this possibility for the stability data $\left(\E,\set{\Pi^k_\eps}_{\eps\in\E}\right)$.
Any linear  order on the set $\E$ must satisfy the conditions
\begin{align}
 &\eps''<\eps'\ \Rightarrow\
\Hom^{\le0}(\Pi^k_{\eps'},\Pi^k_{\eps''})=0,\label{hom condit}\\
&(j',i')<(j'',i'')\  \Leftrightarrow\ (j'+1,i')<(j''+1,i'')\quad \text{for
}(j',i'),\,(j'',i'')\in\E.\label{shift condit}
\end{align}

Identifying the set $\E$ with generators of the semistable subcategories, we see that the first
condition implies
\[\O(k)[i]<\O(k)[j]\quad\text{and }\O(k+1)[i]<\O(k+1)[j]\quad \text{whenever } i<j.\]
On the other hand, \HN\ \eqref{triangles for P1-2} yields $\O(k)<\O(k+1)[-1]$. Therefore due to
the second condition we have
\[\O(k)[i]<\O(k+1)[i-1] \quad\forall\,i\in\ZZ.\]
Now, using the transitivity   of an order, we see that   linear orders on $\E$, satisfying
\eqref{hom condit} and \eqref{shift condit}, depend  on a parameter $p\in
\NN\cup\{0,+\infty\}=\bar\NN$ and have the form
  \[
    \O(k)[i+p]<\O(k+1)[i-1]<\O(k)[i+p+1]\quad\forall\,i\in\ZZ,\quad \text{if }
    p\in\NN\cup\{0\}, \]
    \[\O(k)[i]<\O(k+1)[j]\quad\forall\,i,\,j\in\ZZ,\qquad\text{if }p=+\infty.
\]
Let us denote the set $\E$ with the order  corresponding to $p\in\bar\NN$  by $\E_p$. Since the
existence of finite \HN\ for each object of $D^b(\Coh \PP^1)$ is obvious, we conclude that
$\left( {\E_p},\{\Pi^k_\eps\}_{\eps\in\E_p}\right)$ are the finest  stability data. We call them
exceptional.

\subsection{The set of all stability data for ${\boldsymbol {D}^{\boldsymbol b}\boldsymbol {(}
{\boldsymbol \Coh}{\boldsymbol \PP}^{\boldsymbol 1}{\boldsymbol )}}$}\label{T-STABILIIES FOR P1}
In this subsection we prove the following classification result.

\begin{Thm}\label{All t-stabilities for P1}
 Any finest t-stability on $D^b(\Coh \PP^1)$ is either a  standard one or an exceptional one,
  thus one of the finest t-stabilities constructed in the previous sections.
\end{Thm}
\begin{proof}{.}
Let $\DST \Phi \vf \Pi$ be a t-stability on $D^b(\Coh \PP^1)$ that incomparable with any
standard one. We show that   $\DST \Phi \vf \Pi$ is weaker then one of the exceptional
stabilities.

Recall that the standard stability data $\DST \M \mu \Pi$  depends on a linear order on $\M$ in
the notations of Subsection  \ref{Standard t-stability}. Suppose every $\M$-semistable object is
$\Phi$-semistable. Since each possible linear order on $\M$ gives one of the standard
t-stabilities, we obtain that $\Phi$ is weaker then $\M$. Therefore if $\Phi$ is not comparable
with $\M$, then there exists a $\M$-semistable object that is not $\Phi$-semistable. By
definition   semistable subcategories are closed under extension. Consequently, without loss of
generality, we can assume such $\M$-semistable object is one of the generators of $\Pi_\mu$:
$\O_x[i]$ or $\O(n)[i]$.

Let $X$ be one of the generators. Since $X$ is not $\Phi$-semistable, it has a nontrivial \HN\
w.r.t. $\Phi$. Consider the first \tr\ of the system
\begin{equation}\label{dest triangle}
    \TRM{F^1X}X{X_{\vf_0}}{}{f}\quad.
\end{equation}
We know that $X_{\vf_0}$ is $\Phi$-semistable and
\begin{equation}\label{dest triangle cond}
\Hom^{\le0}(F^1X,X_{\vf_0})=0.
\end{equation}
The shift of the derived  category  acts on the set of all \tr s and on semistable
subcategories. Besides the condition \eqref{dest triangle cond} is invariant under the shift
 too. Thus, we can assume that $X$ is either $\O_x$ or $\O(n)$.

We need the following lemma.
\begin{Lem}\label{dest triangles for P1}
Let  $X\in D^b(\Coh \PP^1)$ is $O_x$ or $\O(n)$, then there exists $k\in \ZZ$ such that the \tr\
\eqref{dest triangle}, satisfying the condition \eqref{dest triangle cond}, is isomorphic to one
of\,\footnote{We identify each sheaf $X\in \Coh\PP^1$ with 0-complex $X[0]\in D^b(\Coh \PP^1)$.}
the triangles \eqref{triangles for P1-1}--\eqref{triangles for P1-3}.
\end{Lem}
Since the \tr\ \eqref{dest triangle} is the first triangle of \HN, we obtain that $x_0\O(k)[0]$
or $x_0\O(k)[1]$ is $\Phi$-semistable for some natural $x_0$. Using Corollary \ref{direct
summands} and Definition \ref{st-data}, we conclude  that    the objects $\O(k)[i]$ are
$\Phi$-semistable for all $i\in\ZZ$, and  the \tr\ \eqref{dest triangle} should have  the form
\begin{equation}\label{first step}
\TR{y_1\O(k+1)[j-1]}X{x_0\O(k)[j]}\quad,
\end{equation}
where $x_0,\,y_1$ are natural numbers and $j=0,\,1$.

Let us show, that $y_0\O(k+1)[j-1]$ is $\Phi$-semistable as well. Suppose not, then   the \HN\
for $X$ continues:
\[\begin{diagram}[size=1.5em]
    && x_0\O(k)[j]&&&&X_{\vf_1}\\
    &\ruTo&&\rdDashto&&\ruTo&&\rdDashto\\
   X&&\lTo&&y_0\O(k+1)[j-1]&&\lTo&&F^2X
\end{diagram}\]
It follows from the definition of \HN\ and Proposition \ref{hom prop of HN}, that
\begin{align}
  &\Hom^{\le 0}(F^2X, X_{\vf_1})=0,\label{hom 2m 1q}\\
 % &\Hom^{\le 0}(F^2X, x_0\O(k)[j])=0,\label{hom 2m 0q}\\
   &\Hom^{\le 0}(X_{\vf_1}, x_0\O(k)[j])=0.\label{hom 1q 0q}
\end{align}
The condition \eqref{hom 2m 1q} and Lemma \ref{dest triangles for P1} allow us to conclude, that
there exist  an integer $m$ and natural numbers  $x_1$, $y_2$ such that
\[\text{either }X_{\vf_1}=x_1\O(m)[j-1],\quad F^2X=y_2\O(m+1)[j-2],\quad \text{if }k<m-1,\]
\[\text{or }X_{\vf_1}=x_1\O(m)[j],\quad F^2X=y_2\O(m+1)[j-1],\quad \text{if }k>m.\]
In the first case ($k<m-1$)
\[\Hom^{0}(X_{\vf_1}, x_0\O(k)[j])=%\Hom^{0}(x_1\O(m)[j-1], x_0\O(k)[j])=
\Ext^1_{\Coh\PP^1}(x_1\O(m),x_0\O(k))\ne0.\] This contradicts \eqref{hom 1q 0q}. In the second
one ($k>m$)
\[\Hom^{ 0}(X_{\vf_1}, x_0\O(k)[j])=\Hom_{\Coh\PP^1}(x_1\O(m),x_0\O(k))\ne0.\]
This contradicts \eqref{hom 1q 0q} again.

Thus there  cannot be the continuation of the \HN, whence $y_0\O(k+1)[j-1]$ is
\mbox{$\Phi$-semistable}. Hence for any $j\in \ZZ$ the object $\O(k+1)[j]$ is
\mbox{$\Phi$-semistable}.

In such a way, we got  $\Phi$-semistable objects $\O(k)[i]$ and $\O(k+1)[i]$ $(i\in\ZZ)$.
Considering more fine stability data (if necessary), we can assume that $\LA\O(k)[i]\RA$ and
$\LA\O(k+1)[i]\RA$ are $\Phi$-semistable subcategories.  Let us show that each $\Phi$-semistable
object belongs to one of them.

Indeed, any object $X$ of the derived category $D^b(\Coh\PP^1)$ could be decomposed into a
direct sum of $i$-complexes $X_i[i]$ (Prop. \ref{dh=1}), where   $X_i$ is a direct sum of
invertible sheaves and torsion sheaves concentrated each in one point. Each these direct summand
of $X$ has \HN\ w.r.t. $\Phi$ with $\O(k)[i]$ and $\O(k+1)[i]$ as quotients (see
\eqref{triangles for P1-1}--\eqref{triangles for P1-2}). Therefore, according to
Proposition~\ref{atom. and glu. t-filtr}, one can construct \HN\ for any object $X$ with only
objects $a\O(k)[i]$, $b\O(k+1)[i]$ as quotients. Now, it follows from the uniqueness of the \HN\
that only objects $x_0\O(k)[i]$ and $x_1\O(k+1)[i]$ are $\Phi$-semistable. Consequently  the
collection of $\Phi$-semistable subcategories coincides with the such collection for some
exceptional t-stability.

The action of the shift on the $\Phi$-semistable subcategories is uniquely determined. Finally,
any admissible linear order on the collection of $\Phi$-semistable subcategories gives us the
exceptional t-stability. This concludes the proof of the theorem.
\end{proof}

\vspace{-3ex}

\begin{Cor}\label{finest t-stab for P1}
    Stability data on $D^b(\Coh\PP^1)$ is weaker then either the standard or the exceptional
    t-stability.
\end{Cor}

\begin{proof}{ of Lemma \ref{dest triangles for P1}.} We start with an observation.

\begin{Rem}\label{nontrivial maps} Suppose
     in the triangle \eqref{dest triangle}   $X_{\vf_0}=X_1\oplus X_2$ and
    $f =f_1\oplus f_2$, where\\ \mbox{$f_i\in\Hom(X,X_i)$}. If $f_2=0$, then $F^1X=Y\oplus X_2$ and
    $\Hom^{\le 0}(F^1X,X_{\vf_0})\ne0$.
\end{Rem}

Let $X$ from triangle \eqref{dest triangle}  be one of the objects $\O_x[0]$ or $\O(n)[0]$.
Taking into account Proposition~\ref{dh=1} and the previous remark, we can assume that
\begin{align*}
 &X_{\vf_0}=X_0[0]\oplus X_1[1]& \text{with } &X_i\in \Coh\PP^1,\\
 &f=f_0\oplus f_1& \text{with }&f_0\in \Hom_{\Coh\PP^1}(X,X_0),\\
 &&&                              f_1\in \Ext^1_{\Coh\PP^1}(X,X_1),\\
 &f_i=0\ \Leftrightarrow\ X_i=0.
\end{align*}
Then $F^1X$ is isomorphic to $C[-1]$, where $C$ is the cone of the morphism $f$.

Consider the extension corresponding to $f_1\in\Ext^1_{\Coh\PP^1}(X,X_1)$:
\begin{equation}\label{ext X1EX for P1}
  0\To X_1\To E\xTo{d_{-1}} X\To0.
\end{equation}
The object $X_1[1]$ is isomorphic to the complex
\[0\To E\xTo{d_{-1}} X\To0,\]
and the morphism $f_1$ can be represented as the morphism of the complexes
\[
\begin{diagram}[size=1.6em]
    &&0&\rTo&X&\rTo&0\\
    &&\dTo&&\dTo>{\scriptstyle id_X}&&\dTo\\
    0&\rTo&E&\rTo^{\scriptstyle d_{-1}}&X&\rTo&0
\end{diagram}\ .\]
Hence the morphism $f:X\To X_0[0]\oplus X_1[1]$ is following:
\[
\begin{diagram}[size=1.6em]
    &&0&\rTo&X&\rTo&0\\
    &&\dTo&&\dTo>{\scriptstyle id_X\oplus f}&&\dTo\\
    0&\rTo&E&\rTo^{\scriptstyle d_{-1}\oplus 0}&X\oplus X_0&\rTo&0
\end{diagram}\ .\]
The cone of this morphism is the complex
\[0\To X\oplus E\xTo{\de_{-1}}X\oplus X_0\To0,\]
where $\de_{-1}$ one can write as the matrix
\[
\de_{-1}=\left(\begin{array}{cc}
    id_X&d_{-1}\\
    f_0&0
    \end{array}
    \right)\,.
\]
The cone complex  is isomorphic to $(\Ker \de_{-1})[1]\oplus (\Coker \de_{-1})[0]$ in the
derived category.

Consider the commutative diagram in the abelian category $\Coh\PP^1$.
\[\begin{diagram}[size=1.6em]
    &&0&&0&&0\\
    &&\dTo&&\dTo&&\dTo\\
0&\rTo&X_1&\rTo&E&\rTo^{\scriptstyle d_{-1}}&X&\rTo&0\\
    &&\dTo&&\dTo&&\dTo&&\dTo\\
    0&\rTo&\Ker\de_{-1}&\rTo&X\oplus E&\rTo^{\scriptstyle \de_{-1}}&X\oplus X_0&\rTo&\Coker \de_{-1}&\rTo&0\\
    &&\dTo&&\dTo&&\dTo&&\dTo\\
    0&\rTo&\Ker f_0&\rTo&X&\rTo^{\scriptstyle f_0}&X_0&\rTo&\Coker f_0&\rTo&0\\
&&\dTo&&\dTo&&\dTo&&\dTo\\
    &&0&&0&&0&&0\\
\end{diagram}\,.\]
We obtain that $\Coker\de_{-1}=\Coker f_0$ and the condition \eqref{dest triangle cond} has the
form
\[
    \Hom^{\le0}\Bigl((\Ker\de_{-1})[0]\oplus (\Coker f_0)[-1],\,X_0[0]\oplus X_1[1]\Bigr)=0,
\]
whence
\begin{align}
  %&\Hom^{\le0}\Bigl((\Coker f_0)[-1],\,X_0[0]\Bigr)=0,\label{CoNd1}\\
%  &\Hom^{\le0}\Bigl((\Coker f_0)[-1],\,X_1[1]\Bigr)=0,\label{CoNd2}\\
  &\Hom^{\le0}\Bigl((\Ker\de_{-1})[0],\,X_0[0]\Bigr)=0,\label{CoNd3}\\
  &\Hom^{\le0}\Bigl((\Ker\de_{-1})[0],\,X_1[1]\Bigr)=0.\label{CoNd4}
\end{align}

Assume that $X=\O_x$ and $X_0\ne0$. Then $\Ker f_0=0$, i.e. $\Ker \de_{-1}=X_1$. This
contradicts \eqref{CoNd4}. Therefore $X_0=0$ and the triangle is obtained from the exact
sequence \eqref{ext X1EX for P1}. In particular, $C[-1]\!=\!E$. So, the condition \eqref{dest
triangle cond} means  $\Hom^{\le 0}(E,X_1[1])=0$, i.e. \mbox{$\Hom_{\Coh\PP^1}(E,X_1)\!=\!0$}
and $\Ext^1_{\Coh\PP^1}(E,X_1)=0$. Besides, by construction, $\Hom_{\Coh\PP^1}(X_1,E)\ne0$. Now
the proof in the case $X=\O_x$ follows from the remark:

\begin{Rem}\label{except pair}
    If sheaves $E$ and $X_1$ on $\PP^1$ satisfy the conditions
    \[\Hom_{\Coh\PP^1}(E,X_1)=0,\
 \Ext^1_{\Coh\PP^1}(E,X_1)=0,\ \Hom_{\Coh\PP^1}(X_1,E)\ne0,\]
 then there exist an integer $k$ and natural numbers  $e$, $x_1$ such that $E=e\O(k+1)$, $X_1=x_1\O(k)$.
\end{Rem}
\begin{proof}{\ } is left to the reader.\end{proof}

 To continue consider $X=\O(n)$. Suppose at first that neither $X_0\ne0$ nor $X_1\ne0$. Then $X_1$ has an
invertible direct summand $L$, since in the converse case $f_1=0$, i.e. $X_1=0$ (see Rem.
\ref{nontrivial maps}). Hence $\Ker \de_{-1}$ has an invertible direct summand $L'$ as well.

Assuming $\Ker f_0=0$, we get $\Ker\de_{-1}=X_1$, which contradicts \eqref{CoNd4}. Therefore
\mbox{$\Ker f_0\ne0$}, i.e.  $\Im f_0$ is a torsion sheaf, and $X_0$ has a nonzero torsion $T$.
Nevertheless the fact that the space $\Hom_{\Coh\PP^1}(L',T)\ne0$ contradicts \eqref{CoNd3}.
Thus in the case $X=\O(n)$ we have
\[\text{either } X_1=0\qquad\text{or } X_0=0.\]

If $X_1=0$ the triangle \eqref{dest triangle} is obtained from the exact sequence
\[0\To\Ker f_0\To\O(n)\xTo{f_0}X_0\To\Coker f_0\To 0\]
and has the form
\begin{equation}\label{f0}
  \TRM{(\Ker f_0)[0]\oplus(\Coker f_0)[-1]\,.}{\O(n)}{X_0[0]}{}{f_0}
\end{equation}
The condition \eqref{dest triangle cond} gives
\begin{align}
  &\Hom^{\le0}((\Ker f_0)[0],X_0[0])=0,\label{CoN5}\\
  &\Hom^{\le0}((\Coker f_0)[-1],X_0[0])=0.\label{CoN6}
\end{align}
As above, the assumption $\Ker f_0\ne0$ contradicts \eqref{CoN5}. Therefore $\Ker f_0=0$.
Furthermore, the condition \eqref{CoN6} and Remark \ref{nontrivial maps} imply that the triangle
\eqref{f0} is isomorphic to \eqref{triangles for P1-2} for some integer $k>n$.

Finally suppose $X_0=0$. Then the triangle \eqref{dest triangle} is obtained from the exact
sequence
\[0\To X_1\To E\xTo{d_{-1}}\O(n)\To  0\]
and has the form
\begin{equation}\label{f1}
  \TRM{E}{\O(n)}{X_1[1]}{}{f_1}
\end{equation}
with $\Hom^{\le0}(E,X_1[1])=0$. Using Remark \ref{except pair} again, we get that the triangle
\eqref{f1} is isomorphic to \eqref{triangles for P1-1} for an integer $k<n-1$. This completes
the proof of the lemma.\end{proof}

\subsection{Classification of t-structures on ${\boldsymbol {D}^{\boldsymbol b}\boldsymbol {(}
{\boldsymbol \Coh}{\boldsymbol \PP}^{\boldsymbol 1}{\boldsymbol )}}$}\label{CLASS T-STRUCT P1}
In this subsection   we give the list of all bounded t-structures on the bounded derived
category of coherent sheaves on $\PP^1$. Notice that the group $\Aut (D^b(\Coh\PP^1))$ of
autoequivalences of the category acts on t-structures and we write only representatives of
t-structures modulo this action. According to \cite{BO} the group $\Aut (D^b(\Coh\PP^1))$ is
generated by the shift, the automorphisms group of $\PP^1$ and Picard's group $\mathrm
{Pic}\,\PP^1$.

As it is well-known, a t-structure $(D^{\ge0},D^{\le 0})$ is uniquely reconstructed by any its
half. Therefore we shall indicate both parts of a t-structure  only if the reconstruction is not
obvious.

The first series  of t-structures is obtained from the standard t-stability (Subsection
\ref{Standard t-stability}) $\DST \M \mu \Pi$.

We start with the tautological t-structure of the derived category:
\[\A^{\le 0}=\LA (\Coh \PP^1)[j],\ j\ge 0\RA.\]
Its core\footnote{Here and further on  we denote the core of a t-structure $(D^{\ge0},D^{\le
0})$ by $D^0$.} is $\A^0=\Coh\PP^1$. This could be shown on a picture.

\begin{center}\unitlength=1mm
\begin{picture}(120,40)
        \put(0,30){$\ldots\,\O[-1]\,\ldots\,\underbrace{\ldots\,\O(-1)[0]\,\O[0]\,\O(1)[0]\,
        \ldots\, \O_x[0]\,\ldots\,}_{\Coh\PP^1[0]}\ldots\,\O(-1)[1]\,\O[1]\,\ldots$}
        \put(22.5,32){\line(0,1){3}}
        \put(22.5,35){\line(1,0){100}}
        \put(100,36){$\A^{\le0}$}
        \put(88.8,32){\line(0,1){2}}
        \put(88.8,34){\line(-1,0){89}}
        \put(7,36){$\A^{\ge0}$}
    \end{picture}
\end{center}\nopagebreak

\vspace{-20mm}

 The next three t-structures are defined by torsion pairs. Recall the definitions and some
 results from \cite{HA, BV}.

 Let $\A$ be an abelian category. A pair $(\A_1,\A_0)$ of full subcategories in $\A$ is called a
 {\itshape torsion pair\/}, if $\Hom(\A_1,\A_0)=0$ and every object $A\in\A$ fits in a unique
exact sequence
\[\seq{A_1}A{A_0}\]
with $A_i\in\A_i$. A torsion pair $(\A_1,\A_0)$ is {\itshape cotilting\/}, if any object in $\A$
is a quotient of an object in $\A_0$.

The torsion pair $(\A_1,\A_0)$ defines a t-structure on $D^b(\A)$ by
\begin{align*}
 &{^p}{D^b}(\A)^{\le0}=\{A\in D^b(\A)^{\le 0}|\ H^0(A)\in \A_1\},\\
 &{^p}{D^b}(\A)^{\ge0}=\{A\in D^b(\A)^{\ge -1}|\ H^{-1}(A)\in \A_0\}.
\end{align*}
By definition the tilting \ ${}^p\A$ of $\A$ w.r.t. $(\A_1,\A_0)$ is the core of this
t-structure. In the case when the torsion pair $(\A_1,\A_0)$ is cotilting, the derived
categories coincide: $D({}^p\A)=D(\A)$.

\smallskip

For $D^b(\Coh \PP^1)$ we have the following cotilting torsion pairs and t-structures.
\begin{align*}
  &\mathcal{B}_0=\LA \O(n),\ n<0\RA,\quad \mathcal{B}_1=\LA \O(n),\ n\ge0;\ \O_x,\
  x\in\PP^1\RA,\\
  &\mathcal{B}^{\le0}=\LA \O(n)[i],\ n\ge0,\ i\ge0;\ \O_x[i],\ x\in \PP^1,\ i\ge 0;\ \O(n)[j],\ n<0,\
  j\ge 1 \RA,\\
  &\mathcal{B}^0=\LA \O(n)[0],\ n\ge0;\ \O_x[0],\ x\in\PP^1;\ \O(n)[1],\ n<0 \RA.
\end{align*}

\begin{center}\unitlength=1mm
    \begin{picture}(150,35)
        \put(0,30){$\ldots\,\O[-1]\,\ldots\,\underbrace{\ldots\,\O(-1)[0]\,\O[0]\,\O(1)[0]\,
        \ldots\, \O_x[0]\,\ldots\,}_{\Coh\PP^1[0]}\ldots\,\O(-1)[1]\,\O[1]\, \ldots$}
        \put(46,32){\line(0,1){3}}
        \put(46,35){\line(1,0){78}}
        \put(100,36){$\mathcal{B}^{\le0}$}
        \put(111,32){\line(0,1){2}}
        \put(111,34){\line(-1,0){112}}
        \put(7,36){$\mathcal{B}^{\ge0}$}
    \end{picture}
\end{center}

\vspace{-28mm}

\begin{align*}
  &\mathcal{C}_0=\LA \O(n),\ n\in\ZZ\RA,\quad \mathcal{C}_1=\LA \O_x,\
  x\in\PP^1\RA,\\
  &\mathcal{C}^{\le0}=\LA \O_x[i],\ x\in \PP^1,\ i\ge 0;\ \O(n)[j]\ n\in\ZZ,\ j\ge1  \RA,\\
  &\mathcal{C}^0=\LA \O_x[0], \ x\in\PP^1;\ \O(n)[1],\ n\in\ZZ \RA.
\end{align*}

  \begin{center}\unitlength=1mm
     \begin{picture}(122,38)
        \put(0,30){$ \ldots\,\O(-1)[0]\,\O[0]\,
        \ldots\,\underbrace{\ldots\, \O_x[0]\,\ldots\,}_{\text{torsion sheaves}}\ldots\,
        \O(-1)[1]\,\O[1]\,\ldots\,\underbrace{\ldots\, \O_x[1]\,\ldots}_{\text{torsion sheaves}}\ldots$}
        \put(35.2,32){\line(0,1){3}}
        \put(35.2,35){\line(1,0){84}}
        \put(100,36){$\mathcal{C}^{\le0}$}
        \put(93,32){\line(0,1){2}}
        \put(93,34){\line(-1,0){93}}
        \put(7,36){$\mathcal{C}^{\ge0}$}
    \end{picture}
\end{center}

\vspace{-20mm}

The last cotilting torsion pair for $\Coh\PP^1$ depends on an arbitrary nonempty subset
$P\subset\PP^1$.
\begin{align*}
  &\mathcal{D}(P)_0=\LA \O_x,\ x\in P\RA,\quad \mathcal{D}(P)_1=\LA \O_y,\
  y\notin P;\ \O(n),\ n\in\ZZ\RA,\\
  &\mathcal{D}(P)^{\le0}=\LA \O_x[i],\ x\in P,\ i\ge 0;\ X[j],\ X\in\Coh\PP^1,\ j\ge1  \RA,\\
  &\mathcal{D}(P)^0=\LA \O_x[0], \ x\in P;\ \O_y[1], \ y\notin P;\ \O(n)[1],\ n\in\ZZ \RA.
\end{align*}

\begin{center}\unitlength=1mm

    \begin{picture}(120,40)
        \put(0,30){$ \ldots\,
        \underbrace{\ldots\, \O_y[0]\,\ldots\,}_{y\in\PP^1\setminus P}\,
        \underbrace{\ldots\,\O_x[0]\,\ldots\,}_{x\in P}
        \,\ldots\,\O[1]\,\ldots\,\underbrace{\ldots\, \O_y[1]\,\ldots\,}_{y\in\PP^1\setminus P}\,
        \underbrace{\ldots\,\O_x[1]\,\ldots\,}_{x\in P}\ldots$}
         \put(27.5,32){\line(0,1){3}}
         \put(27.5,35){\line(1,0){92}}
         \put(100,36){$\mathcal{D}(P)^{\le0}$}
%        %
        \put(92,32){\line(0,1){2}}
        \put(92,34){\line(-1,0){93}}
        \put(7,36){$\mathcal{D}(P)^{\ge0}$}
    \end{picture}
\end{center}

\vspace{-20mm}

It is obvious that t-structures  $\mathcal{D}(P)\sim \mathcal{D}(P')$ modulo $\Aut
D^b(\Coh\PP^1)$ iff $P'=\vf(P)$ for some  $\vf\in\Aut\PP^1$.

\smallskip

The following t-structures are induced by the exceptional stability data
$\left(\E_p,\set{\Pi^0_\eps}_{\eps\in\E_p}\right)$ (see Subsection \ref{Exceptional
t-stability}), where ordering, as we already know, depends on $p\in\bar\NN$. They are sorted
into two  groups: bounded ($\mathcal{E}(p)$ and $\mathcal{F}(p)$) and unbounded ($\mathcal{G}$,
$\mathcal{H}$, and $\mathcal{I}$), and defined as follows.

%\smallskip

\begin{align*}
  &\mathcal{E}(p)^{\le0}=\LA\O[i],\ i\ge p;\ \O(1)[j],\ j\ge-2\RA,\\
  &\mathcal{E}(p)^{ 0}=\LA\O(p)[0],\ \O(1)[-2]\RA,
\end{align*} where $p$ is a nonnegative integer.

\begin{center}\unitlength=1mm
    \begin{picture}(93,38)
        \put(0,30){$\ldots\,\O[p-1]\,\O(1)[-2]\,\O[p]\,\O(1)[-1]\,\O[p+1]\,\O(1)[0]\,\ldots$}
         \put(20,32){\line(0,1){3}}
         \put(20,35){\line(1,0){73}}
         \put(60,36){$\mathcal{E}(p)^{\le0}$}
%        %
        \put(43.6,27){\line(0,1){3}}
        \put(43.6,27){\line(-1,0){44}}
        \put(7,22){$\mathcal{E}(p)^{\ge0}$}
    \end{picture}
\end{center}

\vspace{-20mm}

\begin{align*}
  &\mathcal{F}(p)^{\le0}=\LA\O[i],\ i\ge p;\ \O(1)[j],\ j\ge-1\RA,\\
  &\mathcal{F}(p)^{ 0}=\LA\O(p)[0],\ \O(1)[-1]\RA,
\end{align*}
where $p$ is a nonnegative integer.

\begin{center}\unitlength=1mm
    \begin{picture}(93,40)
        \put(0,30){$\ldots\,\O[p-1]\,\O(1)[-2]\,\O[p]\,\O(1)[-1]\,\O[p+1]\,\O(1)[0]\,\ldots$}
         \put(35.5,32){\line(0,1){3}}
         \put(35.5,35){\line(1,0){57}}
         \put(80,36){$\mathcal{F}(p)^{\le0}$}
%        %
        \put(59,27){\line(0,1){3}}
        \put(59,27){\line(-1,0){60}}
        \put(7,22){$\mathcal{F}(p)^{\ge0}$}
    \end{picture}
\end{center}

\vspace{-20mm}

\begin{align*}
  &\mathcal{G}^{\le0}=\LA \O[i],\ i\ge 0;\  \O(1)[j],\ j\in\ZZ \RA,\\
  &\mathcal{G}^{\ge0}=\LA \O[i],\ i\le 0 \RA,\\
  &\mathcal{G}^{ 0}  =\LA \O[0]\RA.
\end{align*}

\begin{center}\unitlength=1mm
     \begin{picture}(92,40)
        \put(0,30){$\underbrace{\ldots\,\O[-1]\,\O[0]\,\O[1]\,\ldots}_{\O[j],\ j\in\ZZ}\,
        \underbrace{\ldots\O(1)[-1]\,\O(1)[0]\,\O(1)[1]\,\ldots}_{\O(1)[j],\ j\in\ZZ}$}
         \put(16.4,32){\line(0,1){3}}
         \put(16.4,35){\line(1,0){75}}
         \put(80,36){$\mathcal{G}^{\le0}$}
%%        %
        \put(23.7,31){\line(0,1){3}}
        \put(23.7,34){\line(-1,0){24}}
        \put(7,36){$\mathcal{G}^{\ge0}$}
    \end{picture}
\end{center}

\vspace{-20mm}

\begin{align*}
  &\mathcal{H}^{\le0}=\LA  \O(1)[j],\ j\ge 0 \RA,\\
  &\mathcal{H}^{\ge0}=\LA \O [i],\ i\in\ZZ;\ \O (1)[j]\ ,\ j\le 0 \RA,\\
  &\mathcal{H}^{ 0}  =\LA \O (1)[0]\RA.
\end{align*}

\begin{center}\unitlength=1mm
    \begin{picture}(96,40)
        \put(0,30){$\underbrace{\ldots\,\O[-1]\,\O[0]\,\O[1]\,\ldots}_{\O[j],\ j\in\ZZ}\,
        \underbrace{\ldots\O(1)[-1]\,\O(1)[0]\,\O(1)[1]\,\ldots}_{\O(1)[j],\ j\in\ZZ}$}
         \put(59,32){\line(0,1){3}}
         \put(59,35){\line(1,0){37}}
         \put(90,36){$\mathcal{H}^{\le0}$}
%%        %
        \put(72,31){\line(0,1){3}}
        \put(72,34){\line(-1,0){73}}
        \put(7,36){$\mathcal{H}^{\ge0}$}
    \end{picture}
\end{center}

\vspace{-20mm}

\begin{align*}
  &\mathcal{I}^{\le0}=\LA  \O(1)[j],\ j\in\ZZ \RA,\\
  &\mathcal{I}^{\ge0}=\LA \O[i],\ i\in\ZZ \RA,\\
  &\mathcal{I}^{0}  =\LA 0\RA.
\end{align*}

\begin{center}\unitlength=1mm

    \begin{picture}(94,40)
        \put(0,30){$\underbrace{\ldots\,\O[-1]\,\O[0]\,\O[1]\,\ldots}_{\O[j],\ j\in\ZZ}\,
        \underbrace{\ldots\O(1)[-1]\,\O(1)[0]\,\O(1)[1]\,\ldots}_{\O(1)[j],\ j\in\ZZ}$}
         \put(38,32){\line(0,1){3}}
         \put(48,35){\line(1,0){45}}
         \put(90,36){$\mathcal{I}^{\le0}$}
%%        %
        %\put(88,32){\line(0,1){3}}
        \put(48,35){\line(-1,0){50}}
        \put(7,36){$\mathcal{I}^{\ge0}$}
    \end{picture}
\end{center}

\vspace{-20mm}

\begin{Rem}
    It can be easily checked that the cores of the t-structures induced by the exceptional
    stability data are following:
\begin{itemize}
  \item $\mathcal{E}(p)^0$ $\forall p\ge0$ and $\mathcal{F}(p)^0$  $\forall p\ge1$ are equivalent to
   direct sum of the categories of vector spaces $\Vect\oplus\Vect$;
  \item $\mathcal{F}(0)^0$ is equivalent to the category of the quiver
  \begin{picture}(32,10)(-2,-2.5)
    \put(0,0){\circle*{2}}
    \put(27,0){\circle*{2}}
    \qbezier(2,2)(13.5,8)(23,3)
    \put(23,3){\vector(2,-1){2}}
    \qbezier(2,-2)(13.5,-8)(23,-3)
    \put(23,-3){\vector(2,1){2}}
\end{picture}  representations;
  \item $\mathcal{G}^0$ and  $\mathcal{H}^0$ are equivalent to $\Vect$; and
  \item $\mathcal{I}^0$ is the zero category.
\end{itemize}
Since the category $\Vect$ is splitting,    in the case of t-structures induced by the
exceptional stability data  only for $\A=\mathcal{F}(0)^0$ the derived category $D^b(\A)$ is
equivalent to $D^b(\Coh \PP^1)$ (see \cite{BOR}).
\end{Rem}

Now we recall, that a given bounded t-structure on $D^b(\Coh\PP^1)$ induces weak  stability data
(Lemma \ref{stability by t-str}). By Corollary \ref{finest t-stab for P1} this t-stability is
weaker  either the standard ($\M$) one or the exceptional ($\E$) one. Therefore according to
Proposition \ref{class t-struct} any  given bounded t-structure   can obtained by subdividing
  the slope  set $\M$ or $\E$  in
 two disjoint parts. Thereby we proved the following theorem.

 \begin{Thm}\label{list t-structures for P1} The defined above t-structures
 $\A$ --- $\mathcal{F}(p)$ are the only bounded t-structures on $D^b(\Coh\PP^1)$ modulo the group of
 autoequivalences of the category.
 \end{Thm}

\section{Stability data for a  curve of positive genus}\label{ELLIPTIC}

\noindent  In this section we study stability data   for the bounded derived category of
coherent sheaves on a smooth   curve $C_g$ of a  positive genus $g$. We denote the triangulated
category $D^b(\Coh C_g)$ by~$\T_g$.

At first we define the set of standard stability structures extending the Mumford--Takemoto
stability. Then we show that a t-stability for $\T_g$ is induced by a stability on the abelian
category $\Coh C_g$. In the case of elliptic curve $C_1$ we prove that there exists a unique
type of the finest t-stability on $\T_1$. Namely, the standard t-stability depending on ordering
$\bar\QQ$ copies of the curve. Finally we list all bounded t-structures on $\T_1$.

 It was   shown   in  Proposition
\ref{t-stab from core}  that a t-stability on a derived category of an abelian category  $\A$
  can be constructed as an extension of a stability on $\A$. For the abelian category $\A$
   we may define stability structure
 via a positive base of $K_0^*(\A)$  (Definition \ref{abstract slope}). In
the case $\A=\Coh C_g$ the base $(\rk,\deg)$ ($\rk$ is rank and $\deg$ is degree of a sheaf) is
positive. So we have a slope $\gamma(F)$ for every sheaf $F\in \Coh C_g$, where
\[
\gamma(F)=\begin{cases}
   \frac 1 \pi \arcctg\left(-\frac {\deg F}{\rk F}\right)& \text{when } \rk F\ne0,\\
   1& \text{otherwise}.
\end{cases}\]
Since $\arcctg x$ is a strictly decreasing function, this slope is equivalent to
\[
\bar\mu (F)=
  \begin{cases}
    \frac{\deg F}{\rk F} & \text{when }\rk F\ne 0, \\
    +\infty & \text{otherwise}.
  \end{cases}\]
Recall that in these terms a sheaf $E$ is $\bar\mu$-semistable, if for any $0\ne F\subset E$ the
inequality $\bar\mu(F)\le \bar\mu(E)$ is valid. Therefore, any torsion sheaf is
$\bar\mu$-semistable, and a torsion free sheaf on $C_g$ is  $\bar\mu$-semistable if and only if
it is MT-semistable (for Mumford--Takemoto stability).

We shall use the notation  $\bar\QQ$ for the slope set $ \QQ\sqcup\{+\infty\}$ with the natural
order. Thus we obtain the stability  $\DST {\bar\QQ} {q} \Pi$ for the abelian category $\Coh
C_g$ and, furthermore, a t-stability for the derived category $\T_g$. It is clear that the
constructed t-stability is not the finest.

Let us refine $\bar\QQ$. Note that for every rational number $q=\frac d r$ with coprime $d\ne0$
and $r>0$ (or $d=0$, $r=1$) there exists a bundle on $C_g$ of rank $r$ and degree $d$ that is
stable w.r.t. Mumford--Takemoto stability. Denote the set of all such MT-stable bundles by
$\m_q$. It is well-known, that for $F,G\in\m_g$ \quad  $\Hom(F,G)=\Hom(G,F)=0$ whenever $F\ne
G$. If in addition we denote by $\m_\infty$ the curve $C_g$ parameterizing torsion sheaves
$\O_x$, we obtain the set $\M_g$ of semistable subcategories $\Pi_{(q,F)}=\LA F\RA$, where
$q\in\bar\QQ$, and   $F\in\m_q$. Choosing a linear order on each $\m_q$ in arbitrary way and
extending lexicographically  the orders to an order on $\M_g$, we construct the finest stability
$\DST {\M_g} \mu \Pi$ on the abelian category $\Coh C_g$ (because the existence of
Harder--Narasimhan filtration is obvious).

Furthermore, introducing as the slope set $D(\M_g)=\ZZ\x\M_g$ with its lexicographic order and
defining as the semistable subcategories $\Pi_{(i,q,F)}=\LA F[i]\RA$, where $F\in\m_q$, we get
the finest stability data on $\T_g$ (Proposition
 \ref{t-stab from core} and
Proposition \ref{comp. and fine. cond}). We shall call  various t-stabilities obtained this way
the standard t-stabilities.

Note that in the case of an elliptic curve $C_1$   any set $\m_{\frac r d}$ of MT-semistable
sheaves  with coprime $r$ and $d$ consists of indecomposable simple  bundles and   is naturally
isomorphic to $C_1$ (\cite{A}). Therefore the slope set $\M_1$ of a standard stability on $\Coh
C_1$ is the direct product $\bar\QQ\x C_1$.

Further, we prove the following proposition.

\begin{Prop}\label{Induced t-stability}
    Any a t-stability on $D^b(\Coh C_g)$ for $g>0$ is induced by a stability on $\Coh C_g$.
\end{Prop}

\smallskip {\bf Proof.}
 Since homological dimension of $\Coh C_g$ is 1, the proposition is true iff any sheaf $E\in
 \Coh C_g$ considered as a 0-complex has a \HN\ (w.r.t. the t-stability)
 $E\rightsquigarrow(E_{\vf_0},\ldots, E_{\vf_n})$ such that $E_{\vf_i}\in\Coh C_g$ $\forall\,i$.

 Let us consider a \tr
\begin{equation}\label{dest tr for Cg}
  \TRM YEX{}{f}
\end{equation}
with
\begin{equation}\label{dest cond for Cg}
  \Hom^{\le0}(Y,X)=0.
\end{equation}
To verify that all quotients of \HN\ for a sheaf are sheaves as well, it is sufficient to show
that the objects $X$ and $Y$ from the triangle are sheaves, whenever $E$ is.

Now the proposition immediately  follows from the following  lemma.

\begin{Lem}\label{reduction to triple}
    Suppose  $E\in \Coh C_g$ is included in the triangle \eqref{dest tr for Cg}, satisfying the
    condition \eqref{dest cond for Cg}. Then $X,Y\in \Coh C_g$ and $\Hom_{\Coh C_g}(Y,X)=0$.
\end{Lem}
\begin{proof}{ of the lemma.} Since homological dimension of $\T_g$ equals 1, taking into
account Remark \ref{nontrivial maps} we can assume that
 $X=X_0[0]\oplus X_{ 1}[1]$, where $X_0,\,X_{ 1}\in \Coh C_g$, and $f=f_0\oplus
f_{ 1}$ with $f_0\in \Hom_{\Coh C_g}(E,X_0)$, $f_{ 1}\in \Ext^1_{\Coh C_g}(E,X_{ 1})$. Moreover,
the condition $f_i$=0 is equivalent to $X_i=0$.

Consider the extension
\begin{equation}\label{el,extension}
  0\To X_{ 1}\To G\xTo{d_{-1}}E\To 0,
\end{equation}
corresponding to $f_{ 1}$, and realize the object $X_{ 1}[1]$ as a complex $ G\xTo{d_{-1}}E $.
Then $f_{ 1}$  can be represented as the morphism of complexes:
\[
\begin{array}{ccc}
0&\To&E\\
\downarrow&&\downarrow\rlap{$\scriptstyle id_E$}\\
G&\xTo{d_{-1}}&E
\end{array}\ \ \ .\]
Hence the direct sum of maps  $f_0\oplus f_{-1}$ is the morphism:
\[
\begin{array}{ccc}
0&\To&E\\
\downarrow&&\downarrow\rlap{$\scriptstyle id_E\oplus f_0$}\\
G&\xTo{d_{-1}}&E\oplus X_0
\end{array}.\]
A cone  $C(f_0\oplus f_{ 1})$ becomes equal to a complex  of the form
\[ E\oplus G\xTo{\de_{-1}} E\oplus X_0,\]
where   $\de_{-1}$ is determined by the matrix
\[\left(\begin{array}{cc}
id_E&d_{-1}\\
f_0&0
\end{array}\right).\]
Now it follows from the commutative diagram
\[
        \begin{array}{ccccccccccc}
    &&0&&0&&0\\
    &&\downarrow&&\downarrow&&\downarrow\\
 0&\To&X_{ 1}&\To&G&\xTo{d_{-1}}&E&\To&0\\
    &&\downarrow&&\downarrow&&\downarrow&&\downarrow\\
0&\To&H^{-1}&\To&E\oplus G&\xTo{\de_{-1}}&E\oplus X_0&\To&H^0&\To&0\\
&&\downarrow&&\downarrow&&\downarrow&&\downarrow\\
0&\To&\Ker f_0&\To&E&\xTo{f_0}& X_0&\To&\Coker f_0&\To&0\\
&&\downarrow&&\downarrow&&\downarrow&&\downarrow\\
&&0&&0&&0&&0
        \end{array}
\]
that the cone $C(f_0\oplus f_{_1})$ is isomorphic to $H^{-1}[1]\oplus \Coker f_0[0]$. In this
notation the condition  \eqref{dest cond for Cg} can be rewritten  as
\[\Hom^{\le 0}(H^{-1}[0]\oplus \Coker f_0[-1],X_0[0]\oplus X_{ 1}[1])=0.\]
In particular,
\begin{align}
    & \Hom^0(\Coker f_0[-1],X_0[0])=\Ext^1_{\Coh C_g}(\Coker f_0 ,X_0)=0,\label{el,ext1}\\
    & \Hom^0(H^{-1}[0],X_{1}[1])=\Ext^1_{\Coh C_g}(H^{-1} ,X_{ 1})=0.\label{el,ext2}
\end{align}

By definition of a cokernel we have   that $\Hom_{\Coh C_g}(X_0,\Coker f_0)\ne 0$ unless $\Coker
f_0= 0$. On the other hand, the canonical divisor  $K_g$ on a curve of positive genus is either
zero ($g=1$) or effective ($g>1$). Therefore, $\Hom_{\Coh C_g}(X_0,\Coker f_0(K_g))\ne 0$ also.
Whence  by Serre duality theorem we get  that  $\Ext^1_{\Coh C_g}(\Coker f_0, X_0)\ne 0$. This
contradicts the condition    \eqref{el,ext1}, and we conclude that  $\Coker f_0=0$.

We see from the commutative diagram that in the case    $X_{ 1}\ne 0$ the space  $\Hom_{\Coh
C_g}(X_{ 1},H^{-1})$ is not trivial, consequently  $\Hom_{\Coh C_g}(X_{ 1},H^{-1}(K_g))\ne0$.
Using Serre duality theorem again we obtain $\Ext^1_{\Coh C_g}(H^{-1} ,X_{ 1})\ne0$, which
contradicts the condition \eqref{el,ext2}. Thus $X_{ 1}=0$ and the triangle  \eqref{dest tr for
Cg} reduces to
\[\TRM{\Ker f_0}{E}{X_0}{}{f_0}\quad \ .\]
This completes the proofs of the lemma and the proposition.\end{proof}

As a corollary of the proposition we get the following theorem.

\begin{Thm} Any  stability data $\DST \Phi \vf \Pi$ on the bounded
derived category  $\T_1$ of coherent sheaves on a smooth elliptic curve $C_1$ one can refine to
one of the standard stabilities $\DST{D(\M_1)}\mu \Pi$. In particular, any finest t-stability on
$\T_1$ is standard.
\end{Thm}
\begin{proof}{.} Since any stability data on $\T_1$ is induced by a stability on $\Coh C_1$, it is
sufficient to show, that each indecomposable sheaf $E$ (a generator of a standard semistable
subcategory) is $\Phi$-semistable. Suppose that it is not true. Then by Lemma \ref{reduction to
triple} the sheaf $E$ is included in an exact sequence
\[\seq YEX\]
with  $\Hom_{\Coh C_g}(Y,X)=0$. It follows from Serre duality theorem that $\Ext^1_{\Coh
C_g}(X,Y)=0$. Therefore the exact sequence  splits and  $E=X\oplus Y$.
\end{proof}

Finally, using various admissible   decompositions of the slope set $\M_1$, we list all bounded
t-structures on $\T_1$ modulo the group $\Aut\T_1$. Each of them is determined by a cotilting
torsion pair. Therefore we shall describe  only the pairs.

Let $X$ denotes an indecomposable simple sheaf on $C_1$. Denote by $I$ the union of rational
numbers from the segment $[0,1)$ with $\{\infty\}$. Then each $q\in I$ and a subset $P\subset
\m_q\simeq C_1$ of the moduli space of   slope $q$ MT-stable sheaves on $C_1$ (may be empty)
give the cotilting torsion pair $\Bigl(\A(q,P)_1,\,\A(q,P)_0\Bigr)$, where
\begin{align*}
  &\A(q,P)_0=\LA X|\ \bar\mu(X)<q\ \text{or } \bar\mu(X)=q\ \text{and } X\in P\RA,\\
  &\A(q,P)_1=\LA X|\ \bar\mu(X)>q\ \text{or } \bar\mu(X)=q\ \text{and } X\notin P\RA.
\end{align*}
Note that the standard t-structure is obtained when $q=0$ and $P=\varnothing$.

Thus we have the theorem.

\begin{Thm}
    Any bounded t-structure on the bounded derived category $D^b(\Coh C_1)$ of coherent sheaves
    on a smooth elliptic curve (modulo $\Aut ^b(\Coh C_1)$)  is determined by one of the  cotilting torsion pairs
    $\Bigl(\A(q,P)_1,\,\A(q,P)_0\Bigr)$  described above.
\end{Thm}


\begin{thebibliography}{XXXX}
\bibitem{A} Atiyah, M.F. {\itshape Vector Bundles Over an Elliptic Curve\/}. Proc. Lond. Math.
Soc, {\bf VII} (1957) 414-452.
\bibitem{BEI} Beilinson, A.A. {\itshape Coherent Sheaves on $\PP^n$ and Problems in Linear
Algebra\/}. Funk. An., {\bf 12} (1978), p. 68--69.
\bibitem{BOR} A.I.\,Bondal. {\itshape Helices, Representations of Quivers and Koszul
Algebras\/}. Helices and Vector Bundles: Seminare Rudakov. London Math. Soc., Lect. N. Ser. {\bf
148} (1990). p. 75--95.
\bibitem{BO} A.\,Bondal, D.\,Orlov. {\itshape Reconstruction of a variety from the
derived category and groups of autoequivalences\/}.  % (English. English summary)
Compositio Math. {\bf 125} (2001), no. 3, 327--344.


\bibitem{BV} A.\,Bondal and M.\,Van den Bergh. {\itshape Generators and representability of
functors in commutative and noncommutative geometry\/}. Moscow Mathematical Journal v.\,3 (2003)
N\,1. (It is available also in arXiv:math.AG/0204218 v2.)
\bibitem{BR1} Tom Bridgeland. {\itshape Stability conditions on
triangulated categories\/}.\\ arXiv:math.AG/ 0212237.
\bibitem{BR2} Tom Bridgeland. {\itshape Stability conditions on K3 surfaces\/}. arXiv:math.AG/0307164 v1.
\bibitem{GM} S.\,I.~Gelfand and Yu.\,I.~Manin: {\itshape Methods of Homological Algebra\/},
Springer-Verlag (1996).
\bibitem{D} Michael R. Douglas. {\itshape Dirichlet branes, homological mirror
 symmetry, and stability\/}. arXiv:math.AG/0207021.
\bibitem{GIE} Gieseker, D.: {\itshape On the moduli of vector bundles on an algebraic
surface\/}.  Ann. of Math. 106, 45 (1977).
\bibitem{GK} A.\,Gorodentsev, S.\,Kuleshov: {\itshape Helix theory\/}. Preprint MPI 2001 (97).
\bibitem{HA} D.\,Happel, I.\, Reiten, and S.\, Smalo. {\itshape Tilting in abelian categories
and quasitilted algebra\/}, Memoirs of the AMS, vol. 575, Amer. Math. Soc., 1996.
\bibitem{MA1} Maruyama, M.: {\itshape Moduli of stable sheaves I\/}. J. Math. Kyoto 17, 91
(1977).
\bibitem{MA2} Maruyama, M.: {\itshape Moduli of stable sheaves II\/}. J. Math. Kyoto 18, 557
(1978).
\bibitem{RU} A.~Rudakov. {\itshape Stability for an abelian category\/}. J. Algebra 197 (1997),
 no. 1, 231--245
\end{thebibliography}
\end{document}